\documentclass[11pt,a4paper]{article}
\usepackage{beamerarticle}

\mode<article>{\usepackage{fullpage}}
\mode<presentation>{\usetheme{Berlin}}

\usepackage{amsfonts, amsmath}
\usepackage{graphicx}
\usepackage{nccmath}

\usepackage{color}

\usepackage{kpfonts}
\usepackage[cal=boondoxo]{mathalfa}
\input xypic.sty

\usepackage[english]{babel}
\usepackage{pgf}
\usepackage[utf8]{inputenc}
\usepackage{fontenc}

\newcommand{\DMCF}{$\psi$MCF }

\newcommand{\com}[1]{\opt{draft}{\textcolor{red}{
$\LHD$ #1 $\RHD$\marginpar{\textcolor{red}{$\begin{lema}acksquare$}}}}}

\newcommand{\comb}[1]{\opt{draft}{\textcolor{blue}{
$\LHD$ #1 $\RHD$\marginpar{\textcolor{blue}{$\begin{lema}acksquare$}}}}}

\def\qed{\hfill {\large ${\sqcup\!\!\!\!\!\!\sqcap}$}}

\newenvironment{demo}{{\bf Proof }}
{\qed \\}

\newcommand{\re}{\mathbb R}
\newcommand{\ene}{\mathbb N}

\def\ha{{\widehat{\alpha}}}

\newcommand{\<}{\left<}
\renewcommand{\(}{\left(}
\newcommand{\lb}{\label}
\newcommand{\nn}{\nonumber}

\newcommand{\fracc}{\displaystyle\frac}
\newcommand{\ds}{\displaystyle}
\renewcommand{\>}{\right>}
\renewcommand{\)}{\right)}
\newcommand{\flecha}{\longrightarrow}

\newcommand{\eps}{\ensuremath{\varepsilon}}
\def\a{\alpha}
\def\p{\varphi}
\def\r{\rho}

\newcommand{\Ric}{\operatorname{Ric}}

\def\parcial#1#2{\frac{\partial #1}{\partial#2}}

\def\deri#1#2{\frac{d #1}{d#2}}

\newcommand{\uu}{\frak u}

\newcommand{\us}{\underset}

\def\tr{{\rm tr}}
\def\vec{\overrightarrow}
\def\vH{\vec{H}}

\def\parcial#1#2{\fracc{\partial #1}{\partial#2}}

\def\deri#1#2{\fracc{d #1}{d#2}}

\def\bt{\mathbf{t}}

\def\z{{\frak z}}

\def\tr{{\rm tr}}

\newcommand{\bde}{\begin{defi}}
\newcommand{\ede}{\end{defi}}

\newcommand{\be}{\begin{enumerate}}
\newcommand{\ee}{\end{enumerate}}

\newcommand{\ba}{\begin{array}}
\newcommand{\ea}{\end{array}}

\def\Hp{{H_\psi}}
\def\kp{{\kappa_\psi}}

\def\wtau{{\widetilde \tau}}

\def\og{{\overline g}}

\def\oRic{{\overline Ric}}
\def\oM{{\overline M}}
\def\ona{{\overline \nabla}}

\def\oDelta{\overline{\Delta}}

\def\oO{{\overline \Omega}}

\def\wDe{{\widetilde \Delta}}
\def\wg{{\widetilde {\rm g} }}
\def\ws{{\widetilde s }}
\def\wna{{\widetilde \nabla}}
\def\wN{{\widetilde N}}

\def\wwr{{\widetilde r}}

\def\wM{{\widetilde M}}
\def\wg{{\widetilde g}}
\def\wa{{\widetilde \alpha}}

\def\wH{{\widetilde H}}

\def\wF{{\widetilde F}}

\def\wps{{\widetilde \psi}}
\def\wka{{\widetilde \kappa}}

\def\hM{{\widehat M }}
\def\hP{{\widehat P }}

\def\hH{{\widehat H}}

\def\hg{{\widehat {\rm g} }}

\def\hRic{{\widehat Ric}}

\def\hna{{\widehat \nabla}}

\def\hA{{\widehat {A}}}

\def\hna{{\widehat {\nabla}}}

\newtheorem{defi}{\hspace{12pt} Definition}
\newtheorem{teor}{\hspace{12pt} Theorem}
\newtheorem{prop}[teor]{\hspace{12pt} Proposition}

\newtheorem{lema}[teor]{\hspace{12pt} Lemma}
\newtheorem*{lema*}{\hspace{12pt} Lemma}

\newtheorem{nota}{\hspace{12pt} Remark}
\newtheorem{coro}[teor]{\hspace{12pt} Corollary}

\newcommand{\ben}{\begin{enumerate}}
\newcommand{\een}{\end{enumerate}}
\newcommand{\bi}{\begin{itemize}}
\newcommand{\ei}{\end{itemize}}
\newcommand{\bec}{\begin{equation}}
\newcommand{\eec}{\end{equation}}
\newcommand{\beca}{\begin{equation*}}
\newcommand{\eeca}{\end{equation*}}
\newcommand{\bal}{\begin{align}}
\newcommand{\aal}{\end{align}}
\newcommand{\bala}{\begin{align*}}
\newcommand{\aala}{\end{align*}}

\begin{document}

\title{Type I singularities in the curve shortening flow associated to a density}

\author{ Vicente Miquel and Francisco Viñado-Lereu
\thanks{Research partially supported by  the
DGI (Spain) and FEDER  project MTM2013-46961-P. and the Generalitat Valenciana Project  PROMETEOII/2014/064. The second author has been partially supported by a Grant of the Programa Nacional de Formaci\'on de Personal Investigador 2011 Subprograma FPI-MICINN ref: BES-2011-045388 and partially by a contract CPI-15-209 associated to Project PROMETEOII/2014/064}
}

\maketitle

\vspace{-1cm}
\begin{abstract}
We define Type I singularities for the mean curvature flow associated to a density $\psi$ (\DMCF) and describe the blow-up at singular time of these singularities. Special attention is paid to the case where the singularity come from the part of the $\psi$-curvature due to the density. We describe a family of curves whose evolution under \DMCF (in a Riemannian surface of non-negative curvature with a density which is singular at a geodesic of the surface) produces only type I singularities and study the limits of their blow-ups.
\end{abstract}


\section{Introduction }\lb{In} 

\begin{frame}[allowframebreaks]

The mean curvature flow (MCF for short) of an immersion $F_0:M\flecha \oM$ of a hypersurface $M$ in a complete $(n+1)$-dimensional  Riemannian manifold $(\oM,\og)$ looks for a family of immersions $F: M\times I \flecha \oM$ solution of the equation 
\begin{equation}\label{mcf}
\parcial{F}{t} =  \vec{H} = H N, \qquad F(\cdot,0)=F_0
\end{equation}
where $H$ is the mean curvature of the immersion, and we have used the following {\it convention signs for the mean curvature $H$, the Weingarten map $A$ and the second fundamental form ($h$ for the scalar version and $\a$ for its tensorial version)}:

$A X = - \ona_XN$, $\a(X,Y) = \<\ona_X Y, N\> N = \<A X, Y\> N$, $h(X,Y)= \<\a(X Y), N\>$ , for a chosen unit normal vector $N$, and 

$H= \tr A = \sum_{i=1}^n h(E_i,E_i)$, $\vec{H} = \sum_{i=1}^n \a(E_i,E_i) = H\ N$  for a local orthonormal frame $E_1, ..., E_n$  of the submanifold, where $\ona$ denote the Levi-Civita connection on $\oM$. With the same notation $\ona$ we shall indicate the gradient of a function respect to the metric $\og$.

A $(n+1)$-dimensional {\it manifold with density} $(\oM, \og, \psi)$ is a manifold $\oM$ with a metric $\og$ and a function $\psi:\oM\flecha \re$ where, on any $k$-dimensional submanifold $P$ of $\oM$ ($1\le k \le n+1$), we consider the metric $g$ induced by $\og$ but, instead of the canonical volume element $dv_g^k$ associated to the metric $g$, we use the volume element $dv^k_\psi = e^\psi dv^k_g$ induced by the \lq\lq density '' $\psi$. The volume associated to the density $dv_\psi^k$ is called the $\psi$-volume. 
\begin{equation*}
V^k_\psi(P) = \int_P e^\psi \ dv^k_g \equiv  \int_P \ dv^k_\psi. 
\end{equation*}  

Manifolds with density are being actively studied in many contexts. We refer to \cite{mv1} for a short history in the context of mean curvature flow and to the website reference \cite{web}
maintained by Frank Morgan for a huge list of papers dealing with manifolds with densities. 

The natural generalization of the mean curvature  of a hypersurface $M$ in a manifold with density $\oM$ is obtained by the first variation of the $\psi$-volume of $M$. According to \cite{gro}, \cite{mo} and \cite{rocabamo} it is denoted by $H_\psi$ and given (when $\ona\psi$ has sense) by 
$$ H_\psi = H - \<\ona\psi, N\>.$$

When working in the context of a manifold with density, it is then natural to consider mean curvature flows governed by $H_\psi$ instead of $H$. We shall call this flow
\begin{equation}\label{gmcf}
\parcial{F}{t} =  \vec{H_\psi} = H_\psi\ N =  (H - \<\ona\psi, N\>)N,
\end{equation}
the {\it mean curvature flow with density $\psi$} (\DMCF for short).

When $\psi=$constant, we have the standard {\it mean curvature flow} \eqref{mcf}. 

When $\oM$ is a surface and $M$ is a curve, we shall use the notations $\kappa$ and $\kappa_\psi$ for $H$ and $H_\psi$ respectively. In this setting, the mean curvature flow is also called the curve shortening flow.

\end{frame}

Smoczyk, in \cite{sm}, observed that $H_\psi$ is the mean curvature of the warped product  $M\times_{e^\psi} \re$ in $\oM\times_{e^\psi} \re$, which gives as a consequence that: \lq \lq the evolution under the \DMCF \eqref{gmcf} of a hypersurface $M$ of $\oM$ is equivalent to the evolution of the warped product  $M\times_{e^\psi} \re$ in $\oM\times_{e^\psi} \re$ under the MCF \eqref{mcf}". Without explicit reference to densities, this equivalence (with $S^1$ instead $\re$) was used by Angenent (\cite{an89}),  Altschuler, Angenent and Giga (\cite{aag})   and Huisken (\cite{hu90})   to study the behavior of different hypersurfaces of revolution under MCF. Also Smoczyk used a similar approach in  \cite{sm96}.

We start this paper by  writing explicitely and in detail the equivalence (implicit in \cite{sm96})  of many problems related to mean curvature (among them, the MCF) on submanifolds $\hP$ of a riemannian manifold $\hM$ such that there is Riemannian submersion $\pi:\hM\flecha \oM$ and $\hP$ contains all the fibers through points $p\in \hP$,  with the corresponding problems related to the $\psi$-mean curvature (among them the \DMCF) on  submanifolds $P=\pi(\hP)$ in the manifold  $\oM$ with density $\psi$ whose value at every point of $x\in \oM$ is the logvolume of the fiber of the riemannian submersion over $x$. We also show how this equivalence gives new justifications of the definitions of mean and Ricci curvature associated to a density. This is done in section \ref{DS}.

The above remarks are of some help for the main point of this paper: an introduction to the study of Type I singularities for \DMCF on a manifold with density. We begin it in section \ref{BU} by describing the natural definition of type I and the way of doing the blow-up in this context, ending the section with Proposition \ref{bufl}, which states the convergence of a sequence of blow-ups of a type I singular \DMCF to a type I singular \DMCF in the Euclidean space. 

After this, in sections \ref{singI} and \ref{ConBU}, we describe some new situations where the singularities of the \DMCF are of type I and are localized inside the set of singular points of the density $\psi$. 
The setting for these situations is the following:

The ambient manifold $\oM$ will be a complete riemannian surface with a metric $\og$ that can be written as 
\begin{align}\lb{oMog}
\og = dr^2 + e^{2\p(r)} dz^2 \text{ and with sectional curvature }\overline K\ge 0,
\end{align}
where $\p: \re \flecha \re$ is a smooth function satisfying $\p(s)=\p(-s)$, $r$ denotes the $\og$-distance to the curve $r=0$, and $\p(r)\equiv \p\circ r$.
 The existence of such a metric $\og$ over $\oM$ is equivalent to the existence on $\oM$ of a geodesic ($r=0$ in the coordinates where the metric is written) such that the reflection respect to this geodesic $(r,z)\mapsto (-r,z)$ and the reflections  $(r,a-z)\mapsto (r,a+z)$ are isometries. Examples of these surfaces are the ellipsoids of revolution, where the geodesic $r=0$ is an equator (among these examples is the round sphere) and the flat plane.

We consider on $\oM$ a density $\psi$ which depends only on $r$. That means that there is a smooth function, denoted again by $\psi:]0,\infty[\flecha \re$ such that $\psi(x) = \psi(r(x))\equiv \psi(r)(x)$, where the first $\psi$ is the function defined over $\oM$ and the second $\psi$ is the function defined on $]0,\infty[$.  We  shall also demand $\psi$ to satisfy 
\begin{align}\lb{limps3}
&\limsup_{r\to 0}\(\frac{\psi'''}{\psi'} - \frac2{b^2} {\psi'}^2\) \quad \text{ is bounded from above}, \\
&\lim_{r\to 0}\fracc{\psi^{(n)}(r)}{b/r^n}= (-1)^{n-1} (n-1)!, \quad \text{ for some $b>0$ and $n=1,2,3$,}\lb{limps2}
\end{align} 
where $^{(n)}$ denotes the $n$-th derivative respect to $r$. When $b=m\in\ene$, if $\oM\times_{e^{\psi/m}}S^m$ (where $S^m$ is considered with its standard metric of sectional curvature $1$) is a Riemannian manifold without singularities then  $\psi$  satisfies the conditions \eqref{limps3} and \eqref{limps2} for all $n\in\ene$  (then  $\oM\times_{e^{\psi/m}}S^m$ is a rotationally symmetric space  in the sense of \cite{cami2} and \cite{cami3}). This fact and the equivalence between the \DMCF of $M$ in $\oM$ and the MCF of $M\times_{e^{\psi/m}}S^m$ in $\oM\times_{e^{\psi/m}}S^m$   motivate to consider the conditions \eqref{limps3} and \eqref{limps2} also when $b\notin \ene$.  We remark that, when $\psi$ satisfies \eqref{limps2} and $b$ is not a natural number, $\oM\times_{e^{\psi/k}}S^k$ has singularities whatever $k\in \ene$ be. Then the hypothesis on $\psi$ includes many situations (all when $b\in\re^+-\ene$) where the \DMCF seems to be special and not equivalent to the MCF in any regular Riemannian manifold.  More details are given in section \ref{singI}.

For the initial curve $M_0=F(M,0)$ we shall consider two possibilities:

(i) $M=S^1$ and $M_0$ is a simple closed curve in $\oM$

(ii) $M=[b_1,b_2]$, $M_0$ is simple and there is a region $G=\{(r,z)\in \oM;\ a_1\le z \le a_2\}$ such that $M_0$ is contained in $G$, with $\partial M_0= \{F(b_1,0), F(b_2,0)\} \subset \partial G$ and $M_0$ orthogonal to $\partial G$ at the points $F(b_i,0)$.

Next pictures show examples of these cases when $\oM$ is the round sphere $S^2$.

\includegraphics[scale=0.22]{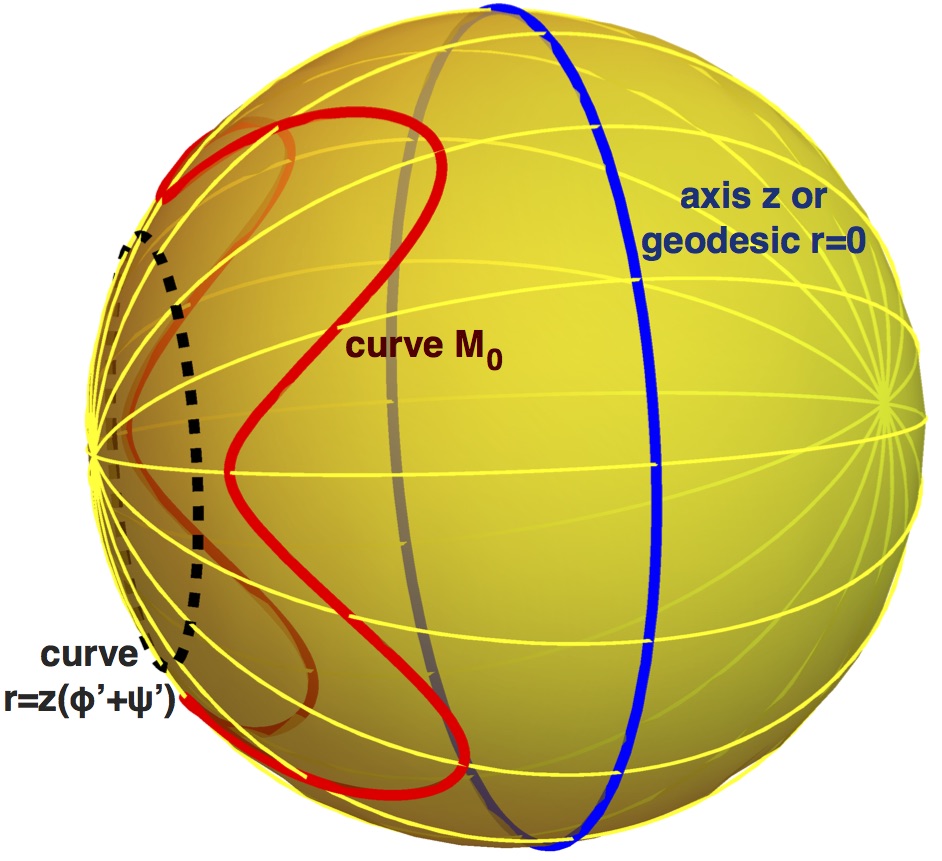} \includegraphics[scale=0.22]{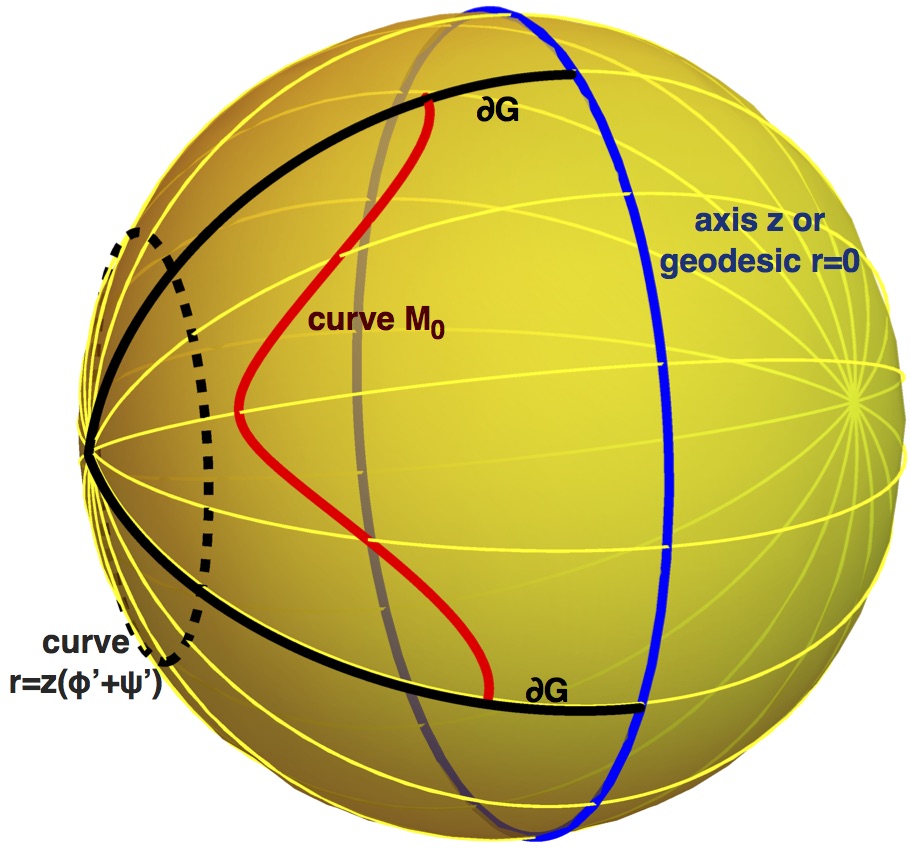}

Along sections \ref{singI} and \ref{ConBU} we shall prove the following

\begin{teor}\lb{thI} Let $(\oM, \og, \psi)$ be a Riemannian surface with density satisfying the conditions \eqref{oMog}, \eqref{limps3} and \eqref{limps2}. Let $M_t$ be the solution of the \DMCF \eqref{gmcf} on a maximal interval $[0,T[$ such that the initial condition $M_0$ is a graph over the geodesic $r=0$, is contained in the band limited by $r=0$ and $r=\min\{\z(\p'+\psi'), \sup\{r; \(\psi''+\psi'^2/b\)|_{[0,r]}\le 0\}\}$, satisfiea  $\kappa_\psi \ge 0$ (but not identically $0$) and (i) or (ii) above in  $(\oM, \og, \psi)$. In case (ii) we add on the problem \eqref{gmcf} the boundary condition: \lq\lq $M_t$ intersects $\partial G$ orthogonally at the boundary of $M_t$ for every $t\in [0,T[$\rq\rq. Then
\begin{enumerate}
\item $\kappa_\psi > 0$ for every $t\in]0,T[$.
\item $M_t$ is a graph over $r=0$ for every $t\in [0,T[$.
\item  $T<\infty$ and the flow $F(\cdot,t)$ is of type I in the sense of Definition \ref{defTId}.
\item At each singular point, a blow-up centered at this point gives a limit flow $\wM_t$ in $\re^2$ with its Euclidean metric and density $\wps^\infty = \ln r^b$ which is a graph over $r=0$ for every time and, after doing a new blow-up,  converges to a $\ln r^b$-shrinker in $\re^2$, which is the line $r={\rm constant}$ in case $b=m\in \ene$.
\end{enumerate}  
\end{teor}
The concept of $\ln r^b$-shrinker which appeared above is
\begin{defi}\lb{fsh}
Given any function $f:\re^{n+1}\flecha \re$, by a $f$-shrinker in $\re^{n+1}$ we understand a hypersurface $F:M\flecha \re^{n+1}$ satisfying $H_f + \<F,N\> =0$.
\end{defi}
\begin{nota} [On the hypotheses $r< \r :=\min\{\z(\p'+\psi'), \sup\{r; \(\psi''+\psi'^2/b\)|_{[0,r]}\le 0\}\}$] 
As it is very usual, $\z(f)$ denotes the first positive zero of a function $f$. 

The inequality $\z(\p'+\psi')>0$ follows from \eqref{limps2}. 

When  $b=m\in \ene$, that is, when the \DMCF is equivalent to a MCF in a Riemannian manifold $\oM\times_{e^{\psi/m}}S^m$, one has that if the sectional curvatures of the planes containing $\partial_r$ and a vector tangent to $S^m$ in $\oM\times_{e^{\psi/m}}S^m$ are non-negative, then always $\psi''+\psi'^2/m\le 0$. 

Moreover the condition $r < \z(\p'+\psi')$ is necessary in order that the lines or circles at distance $\r-\eps$ from the axis $r=0$  collapse to this axis. Again details are given in section \ref{singI}. 
\end{nota}

This remark gives immediately the following Corollary \ref{corThI}, which requires the definition of rotationally symmetric spaces and hypersurface, that we borrow from \cite{cami2} and recall here with  the notation used in this paper.

\begin{defi}\lb{AE}
A rotationally symmetric space  respect to an axis $z$ is a smooth Riemannian manifold $(\hM,\hg)$ admitting cylindrical coordinates $(r,z,u) \in I\times J \times  S^{m}$ respect to which  $\hg$ can be written in the form
\begin{equation}\label{mRSS}
 \hg:= dr^2 + e^{2 \p}(r)\  dz^2 + e^{2\psi/m}(r) \ g_{S},
\end{equation}
where  $g_{S}$ is the standard metric of sectional curvature $1$ on the sphere $S^{m}$.

The curve $r=0$ is a geodesic of $\hM$ and it is called \lq\lq axis $z$ \rq\rq or \lq\lq axis of revolution\rq\rq.
\end{defi}

\begin{defi}\lb{AE2} Let $(\hM,\hg)$ be a rotationally symmetric space respecto to an axis $z$.
A hypersurface of revolution $S$ of $\hM$ generated by a graph over the axis $z$ is a submanifold of $\hM$ that can be described in cylindrical coordinates by the immersion $J\times S^m \flecha \hM / (z,u) \mapsto (r(z),z,u)$, where  $r(z)$ is a smooth function.
\end{defi}

Now, we can state the announced 

\begin{coro}\lb{corThI}
Let $(\hM, \hg)$ be a rotationally symmetric space of dimension $m+2$ and non-negative sectional curvature. Let $M_0$ be a hypersurface of revolution of $\hM$ generated by a graph over the axis \lq\lq z'', with non-negative mean curvature (but not identically $0$)  and contained in the region limited by the cylinders $r=0$ and $r=\z(\p'+\psi')$. Let $M_t$ be the solution of the MCF \eqref{mcf} on a maximal interval $[0,T[$ with the initial condition $M_0$. We consider two cases:

(i) $M_0$ is a closed hypersurface

(ii) $M_0$ is a compact hypersurface with boundary contained in the boundary of a band $G$ limited by two hypersurfaces $z=b_1$ and $z=b_2$. In this case we add to \eqref{mcf} the boundary condition: \lq\lq $M_t$ intersects $\partial G$ orthogonally at the boundary of $M_t$ for every $t\in [0,T[$\rq\rq. 

\noindent Then
\begin{enumerate}
\item $H > 0$ for every $t\in]0,T[$.
\item $M_t$ is a hypersurface of revolution of $\hM$ generated by a graph over the axis \lq\lq z'' for every $t\in [0,T[$.
\item  $T<\infty$ and the flow $F(\cdot,t)$ is of type I.
\item At each singular point, a blow-up centered at this point gives a limit flow which is a graph over $r=0$ for every time and, after doing a new blow-up,  converges to a cylinder in $\re^{m+2}$.\end{enumerate}  
\end{coro}
This Corollary generalizes to manifolds with non-negative curvature some results obtained by G. Huisken (\cite{hu90}) and S. Altlschuler, S. Angenent and Y. Giga (\cite{aag}) for surfaces of revolution in the Euclidean $\re^{n+1}$.

\section{Manifolds with density and Riemannian submersions}\lb{DS}

Let $\pi:\hM \flecha \oM$ be a Riemannian submersion with $m$-dimensional fibers, and $\oM$ of dimension $n+1$. For every $p\in \oM$, the riemannian volume element $dv^{n+1+m}_\hg$ of $\hM$ can be written (with an obvious abuse of notation) as $dv^{n+1+m}_\hg(x)= dv^{n+1}_\og(p) \wedge dv^{m}(x)$ for every $x\in \pi^{-1}(p)$ and $p\in \oM$, where $dv^{n+1}_\og$ denotes the volume element of $\oM$, and $dv^{m}(x)$ denotes the volume element of the fiber $\pi^{-1}(p)$.  As a consequence, for every $\ell$-dimensional submanifold $P\subset\oM$, the volume element of $\hP=\pi^{-1}(P)$ is $dv^{\ell+m}_\hP(x)= dv^{\ell}_g(p) \wedge dv^{m}(x)$, where $dv^\ell_g$ is the volume element of $P$. 

Let us suppose that the fibers $\pi^{-1}(p)$ have finite volume. Define $\psi:\oM\flecha \re$ by 
\bec\lb{defp}
e^\psi(p) := \int_{\pi^{-1}(p)} dv^{m}(x),
\eec
 which defines a density over $\oM$.
One has 
\begin{align}
V^{\ell+m}(\pi^{-1}(P)) &= \int_{\pi^{-1}(P)} dv^{\ell}_g(p) \wedge dv^{m}(x) \nn\\
&= \int_P \(\int_{\pi^{-1}(p)} dv^{m}(x)\) dv^{\ell}_g(p) = \int_P e^\psi dv_g^{\ell} = V_\psi^\ell(P);   
\end{align}
that is, 
{\it the volume of $\pi^{-1}(P)$ coincides with the $\psi$-volume of $P$} .

Now, we relate other geometric invariants associated to the density with geometric invariants in $(\hM,\hg)$.

First, let us compute $\ona\psi$. For every $v\in T_p\oM$, $p\in \oM$, let $c(t)$ be a curve in $\oM$ satisfying $c(0)=p$ and $c'(0)=v$. Let us denote by $\vH_p$ the mean curvature vector of the leaf $\pi^{-1}(p)$ in $\hM$ and by $v^*$ the horizontal lift of $v$ on the fiber $\pi^{-1}(p)$ (that is, $v^*(x) = {\pi^\bot_*}^{-1}v$, where $\pi^\bot_*$ is the restriction of $\pi_*$ to the complementary orthogonal of $T_x\pi^{-1}(p)$). By the definition of the gradient,
\begin{align}\lb{napH}
\<\ona\psi, v\> = \deri{\psi\circ c}{t} \Big|_{t=0} = \deri{}{t}\Big|_{t=0} \ln \int_{\pi^{-1}(c(t))} dv^m = - \frac1{\int_{\pi^{-1}(p)} dv^m} \int_{\pi^{-1}(p)} \<\vH_p(x), v^*\> dv^m(x),
\end{align}
where  the last equality follows from the formula of the first variation of the area applied to the variation $(x,t)\mapsto c_x^*(t)$, $c_x^*(t)$ being the horizontal lift of $c(t)$ starting at $x\in \pi^{-1}(p)$.

Equation \eqref{napH} says that {\it the gradient of $\psi$ at $p\in\oM$ is the averaged mean curvature of the leaf $\pi^{-1}(p)$}.  Moreover, 

\begin{prop}\lb{ponapH2}
If the leaf $\pi^{-1}(p)$ has constant mean curvature, in the sense that its mean curvature vector in $\hM$ is the horizontal lift of a vector tangent to $\oM$, then 
\begin{align}\lb{onapH2}
\ona \psi(p) = - \pi_{*x} \vH_p(x) \text{ for every } x\in\pi^{-1}(p).
\end{align} 
\end{prop}
If $\hna$ and $\ona$ denote, respectively, the covariant derivative in $\hM$ and $\oM$, the equations of a Riemannian submersion give the relation
$\pi_*(\hna_{X^*} Y^*) = \ona_XY$, or, equivalently, $(\hna_{X^*} Y^*)^\bot =( \ona_XY)^*$ for every vector fields $X$, $Y$ defined on $\oM$.

As a consequence, if $M$ is an immersed hypersurface of $\oM$, the second fundamental form $\ha$ of $\pi^{-1}(M)$ in $\hM$ is related with the corresponding one $\a$ of $M$ in $\oM$ by

$\ha(X^*, Y^*) = \a(X,Y)^*$ and $\hH(x) = H(p) + \<\vH_p(x),N^*(x)\>$ for every $x\in\pi^{-1}(p)$.

Then, if $\vH_p$ is the horizontal lift of a vector tangent to $\oM$, the above remark and \eqref{onapH2} give

$\hH(x) = H(p) - \<\ona \psi(p),N(p)\>$.

That is, 
\begin{prop} If $\vH_p$ is the horizontal lift of a vector tangent to $\oM$, the $\psi$-mean curvature of $M$  in $\oM$ coincides with the mean curvature of $\pi^{-1}(M)$ in $\hM$.
\end{prop}

A consequence of these observations is that  the isoperimetric problem for densities, the classification of submanifolds of constant $\psi$-mean curvature and the mean curvature flow for densities on a manifold $\oM$ are equivalent to the corresponding purely Riemannian problems for manifolds $\hM$ such that there is a Riemannian submersion $\pi:\hM\flecha \oM$ with leafs which have \lq\lq constant'' mean curvature vector and submanifolds $S$ of $\hM$ which contain all the fibers in $\oM$ of points $p\in M= \pi(S)$. When a group $G$ acts as a group of isometries on a manifold $\hM$, a $G$-equivariant problem is a problem where only $G$-equivariant domains or $G$-equivariant submanifolds are considered. These $G$-equivariant problems fit in the above class of Riemannian problems. An example is the equivariant MCF studied in \cite{sm96}, which is solved using the equivalent problem on densities, although, once again,  the word density is not mentioned.

\medskip

Let us remark that, not only problems with Riemannian submersions give rise to a problem with densities. When the density function $\psi$ is regular, it is also true the assertion in the opposite direction: {\it any problem with densities is equivalent to many problems on Riemannian submersions}. In fact, given the manifold with a density $(\oM,\og,\psi)$ and any $m$-dimensional Riemannian manifold $Q$ of finite volume $V$, take $\phi=\(e^\psi/V\)^{1/m}$ and the warped product $\hM=\oM\times_{\phi} Q$ defines  a Riemannian submersion $\pi:\hM\flecha \oM; \ \pi(p,x)=p$ for which the formula \eqref{defp} gives the original density $e^\psi$ as the density over $\oM$ defined by $\pi$. 

The hypothesis that fibers have finite $m$-volume may be changed by the hypothesis that the volume forms of the fibers are homothetic, that is: there is a function $e^\psi: \oM \flecha\re$ such that the volume elements of the fibers have the form $dv^m_p = e^\psi(p) {\p_\a}^* d\omega^m$ for every fibered chart $\p_\a$ of a certain atlas of the fibration $\pi:\hM\flecha \oM$ and a fixed volume form  $d\omega^m$ on the canonical fiber.  With this alternative hypothesis,  we still have that $dv^{\ell+m}(x) = dv^\ell(p) \wedge dv^m(x)$, then (up to the isomorphisms $\p_\a^*$) $dv^{\ell+m}(x) = e^\psi(p) dv^\ell(p) \wedge d\omega^m$. This gives the same relations between $H$ and $\hH$ than before. Warped products $\hM = \oM\times _{e^{\psi/m}} Q$ (with any Riemannian manifold $Q$) provide a family of Riemannian submersions $\pi: \hM \flecha \oM$ which satisfy the above alternative hypothesis. When volume$(Q)$ is finite, one obtains for the $\psi$ defining the homothecies the expression $e^\psi(p) = V(\pi^{-1}(p))/V(Q)$, which differs from  \eqref{defp} only in the product by a constant, which does not affect the geometry of the problems.

These considerations give another justification of some of the definitions   of the Ricci curvature in a Riemannian manifold with density. The known formulae for the curvature for a warped product (cf. \cite{ON}) state:
\begin{align}\lb{defRic4}
\hRic(X,Y) &= \oRic(X,Y) - \frac{m}{e^{\psi/m}} \ona^2(e^{\psi/m})(X,Y) \nn \\
&= \oRic(X,Y) - \frac{m}{e^{\psi/m}} \frac{e^{\psi/m}}{m}(\ona^2\psi +\frac1m  \ona\psi\otimes \ona\psi)(X,Y) \nn \\
&= \oRic(X,Y) -(\ona^2\psi + \frac{1}{m} \ona\psi\otimes \ona\psi)(X,Y)
\end{align}  
 which coincides with the Ricci curvature with density in Bayle's thesis (cf. \cite{BayT}), and is called many times the Bakry-Emery tensor. If our  starting subject is the manifold with density $(\oM,\og,\psi)$, the $m\in \ene$ is arbitray , we can take $m\to\infty$ and obtain another usual definition of Ricci curvatue associated to a density (see, for instance, \cite{mo}).

\section{Type I singularities for the  \DMCF and their blow-up}\lb{BU}

 To introduce the concept of Type I flow for the \DMCF, we shall start using the equivalence between \DMCF and certain MCF   stated in the previous section. 

Let $(\oM,\og,\psi)$ be a Riemannian manifold with density. Let $M_0$ be a hypersurface of $\oM$ and $M_t$ be the \DMCF on $\oM$ with $M_0$ as initial condition. We know that this is equivalent to the MCF of $M_0\times_{e^{\psi/m}} Q$ in $\oM\times_{e^{\psi/m}} Q$, and that $\oM\times_{e^{\psi/m}} Q$ could be a singular Riemannian manifold at the points where $\psi$ has singularities. We emphasize that it could happen that $\psi$ be singular and  $\oM\times_{e^{\psi/m}} Q$  be still a regular Riemannian manifold,  as we shall see soon. For the MCF it is known that, if $T$ is the maximal existence time, one has that either the evolution attains a singular point of the metric of $\oM\times_{e^{\psi/m}}Q$ at time $T$ or $|\hA_t|^2$ becomes infinite at $T$ and, then,
\begin{align}\lb{maxAl}
\max_{x\in \hM_t}|\hA_t|^2 \ge \fracc{1}{2(T-t)}.
\end{align}
The second fundamental form $\ha_t$ and the Weingarten map $\hA_t$ of $M_t\times_{e^{\psi/m}} Q$ in $\oM \times_{e^{\psi/m}} Q$ satisfy the equations (see \cite{ON}) 
\begin{align}\lb{hAA}
\hA_t X^*=(A_t X)^*, \quad \hA_t V= - \hna_{V}N_t^* = -\frac{N_t e^{\psi/m}}{e^{\psi/m}} V = -\frac1{m}\<\ona \psi, N_t\>  V
\end{align}
for every vertical vector field $V$. Then   $|\hA_t|^2 = |A_t|^2 + \fracc1{m} \<\ona \psi,N_t\>^2$ and $|\hA_t|^2$ becomes infinite at $T$ if and only if $|A_t|^2$ or $\<\ona \psi,N_t\>^2$ becomes infinite at $T$. This condition means that the hipersurface $M_t$ becomes singular ($\lim_{t\to T}|A_t|^2\to \infty$) or the equation \eqref{gmcf} itself becomes singular when $t\to T$, then it is just the condition for the first singular time of the \DMCF. That is, {\it the final times $T$ for the flows \eqref{mcf}  and \eqref{gmcf} coincide if we do not worry about the possible singularities of the Riemannian manifold $\oM\times_{e^{\psi/m}}Q$}. Of course, if  $\<\ona \psi,N_t\>^2\to \infty$ when $t\to T$, the hipersurface touches the singular points of $\psi$ in the limit when $t\to T$, however, a hypersurface could contain singular points of $\psi$ and keep  $\<\ona \psi,N_t\>^2$ bounded. According to \eqref{maxAl} and \eqref{hAA}, if $T$ is the first singular time of the \DMCF, one has 
\begin{align}\lb{maxAlDG}
\max_{x\in M_t}\(|A_t|^2 + \frac1m\<\ona \psi,N_t\>^2 \)\ge \fracc{1}{2(T-t)}.
\end{align}

The inequality \eqref{maxAlDG} could come from $|A|\to \infty$ or from $|\<\ona \psi,N\>|\to \infty$ or both together. There are interesting  situations where the singularities come from the second case, because this localizes the possible singularities at the singular points of the density $\psi$. As we shall prove in section \ref{singI}, the situations described in Theorem \ref{thI} fall in this type. 

 As in the MCF,  because of property \eqref{maxAlDG}, it has sense to define
 
 \begin{defi}\lb{defTId} A \DMCF is of type I if there is a constant $C>0$ such that 
 \begin{align}\lb{defTI}
&\sup_{M_t}\(|A|^2 + \frac1b\<\ona\psi, N\>^2\) \le \frac{C}{T-t} \qquad \text{ where } b>0 
\end{align}
\end{defi}
Then, if $\ona\psi$ is bounded in the region of $\oM$ where $M_t$ evolves by the \DMCF, the evolution of $M_t$
 is of type I if and only if there is a constant $C>0$ such that $\displaystyle \sup_{M_t}|A|^2 \le \dfrac{C}{T-t}$. 
 
 The number $b>0$ is irrelevant, it only changes the constant $C$ in \eqref{defTI}, but it is the memory of the $m$ in \eqref{maxAlDG} and it is useful in the context of the hypothesis \eqref{limps2} used in Theorem \ref{thI}.

An immediate property of type I evolutions is

\begin{prop}\lb{limF}
If $M$ is compact and $F:M \times [0,T[ \flecha \oM$ is a type I evolution under \DMCF, with $T<\infty$, then $F(\cdot,t)$ converges uniformly, as $t\to T$, to some continuous function $F(\cdot,T): M \flecha \oM$. Moreover, the limit $F(M,T)$ is also compact.
\end{prop}
\begin{demo}
For every $p\in M$, $0\le s < t < T$ one has
\begin{align}
d(F(p,s), F(p,t)) & \le L_s^t(F(p,\sigma)) = \int_s^t \left|\parcial{F}{\sigma}\right| d\sigma = \int_s^t |H-\<\ona\psi,N\>| d\sigma \nn \\
&\le \int_s^t \(|H| + |\<\ona\psi,N\>|\)  d\sigma \le \int_s^t \frac{\sqrt{C}(\sqrt{n}+\sqrt{b})}{\sqrt{T-\sigma}} d\sigma  \nn \\
&= 2 \sqrt{C}(\sqrt{n}+\sqrt{b})\  \(-\sqrt{T-t}+\sqrt{T-s}\) .\lb{dFF}
\end{align}
This shows that the family $F(\cdot,t)$ satisfies Cauchy condition for the topology of the uniform convergence, which proves that the functions $F(\cdot,t)$ converge to a continuous function $F(\cdot,T)$ as $t\to T$, then $F(M,T)$ is compact if $M$ is.
\end{demo}

As in \cite{hu90}, this proposition led us to the 
\begin{defi}
We say that $p\in \oM$ is a blow-up point for the \DMCF of $M$ if there is a $x\in M$ such that $F(x,t)$ converges to $p$ and $|A|(x,t)$ or $|\<\ona\psi,N\>|(x,t)$ become unbounded as $t\to T$. 
\end{defi}

  \subsection{Blow-up of the ambient space}
  
 Let $M_t$ be the evolution of a type I \DMCF of a compact hypersurface $M$ without boundary or with the boundary in the boundary $\partial G$ of a domain $G$ in $\oM$. Let $p\in G\setminus \partial G$ be a blow-up point in $\oM$ and let $x\in M$ such that $F(x,t) \us{t\to T} \to p$. Let us consider a sequence of times   $0<t_1 < t_2 < ... < t_j < ... <T$ which converge to $T$. At each $t_j$, let us rescale the metric $\og$ as
\begin{align}\lb{ogj}
 \og^j := \lambda^2_j \og, \text{ where } \lambda^2_j\equiv \lambda(t_j)^2:=   \fracc{C}{T-t_j}.
 \end{align}

This produces the same rescaling for the metric $g_t^j$ induced on $M_t$ by $\og^j$ respect the metric $g_{t}$ of $M_{t}$ induced by $\og$, that is
$$g_{t}^j = \lambda_j^2\ g_{t}.$$

Let $R>0$ be lower than the injectivity radius of $\oM$ and also lower than the distance from $p$ to $\partial G$.
The pointed Riemannian manifolds $(B^{\og}_{R}(p), \og^j, p)$ (where $\og^0:= \og$) converge, in the Cheeger-Gromov $C^\infty$ topology, to the pointed Euclidean space $(\re^{n+1}, g^e, 0)$. In fact, let us denote by $\Phi_j:(\re^{n+1},g^e)\flecha (T_p\oM, \og^j_p)$ an isometry between these two Euclidean vector spaces; there is a family of open sets $B^e_{\lambda_j R}$ of $\re^{n+1}$ that contain $0$  and a family of diffeomorphisms $\p_j : B^e_{\lambda_j R} \to B^{\og}_{R}(p)$, $\p_j(v)=\exp^{\og^j}_p\circ \Phi_j(v)$ such that $\p_j^*\og^j$  converge to $\og^e$ $C^\infty$-uniformly on every compact of $(\re^{n+1},\og^e)$, because $\p_j^{-1}$ gives the normal coordinates around $p$ of the points in $(B^\og_R(p), \og^j)$ and, computing in these coordinates, we obtain  that the expressions of the metrics $\og^j$ in normal coordinates around $p$ are
\begin{align}
\wg^j_{ik}(0) : = (\p_j^*\og^j)_{ik}(0) &= (\p_0^*\og^0)_{ik}(0) = g^e_{ik}= \delta_{ik}, \\
\wg^j_{ik}(v) : = (\p_j^*\og^j)_{ik}(v) &= (\p_0^*\og^0)_{ik}(v/\lambda_j), \\
\parcial{}{v^\ell}\wg^j_{ik}(v) &= \frac1\lambda_j \parcial{}{v^\ell}\wg^0_{ik}(v/\lambda_j) , \\
\cdots \ & \cdots \ \cdots  \\
\parcial{^r}{v^{\ell_1} \cdots \partial v^{\ell_r}}\wg^j_{ik}(v) &= \frac1{\lambda_j^r} \parcial{^r}{v^{\ell_1} \cdots \partial v^{\ell_r}}\wg^0_{ik}(v/\lambda_j),
\end{align}
which gives the desired smooth convergence of $\p_j^*\og^j$ to $g^e$ on every compact when $j\to\infty$ (equivalent to $\lambda_j\to \infty$).
 
 \subsection{Blow-up of the density function}\lb{bwdf}
 For the gradients of the funtion $\psi$ respect to the metrics $\og$ and $\og^j$, we have 
 \begin{align}
\og(\ona\psi ,X)&= d\psi(X) =  \og^j(\ona^j\psi ,X) = \lambda_j^2 \og(\ona^j\psi ,X),\\
\text{then }   \ona\psi &= \lambda_j^2 \ona^j\psi.
\end{align} 

To consider the function $\psi$ defined on the subsets $B^{e}_{\lambda_j R}$ of $\re^{n+1}$, we have to do the pull-back by $\p^*_j$, and we obtain the induced functions
\begin{align}\lb{defpsij}
\wps^j := \p_j^*\psi = \psi\circ\p_j: B^{e}_{\lambda_j R} \subset \re^{n+1} \flecha \re.
\end{align}
 Let us observe that if $\Phi_0^{-1}(X)=v$, then $\Phi_j^{-1}(X)= \lambda_j v$,  which gives  
 \begin{align}\lb{wpsj}
 \wps^j(v) = \psi(\exp_p^{\og^j}\Phi_j(v)) = \psi(\exp_p^{\og}\Phi_0(v/\lambda_j)) =\wps^0(v/\lambda_j).
 \end{align}  

If $\psi$ is a well defined function on all $\oM$ (without singularities), $\wps^\infty(v):=\lim_{j\to \infty} \wps^j(v)= \lim_{j\to \infty} \wps^0(v/\lambda_j) = \wps^0(0) = \psi(p)$ is a constant function. Here, as above, the limits are in the topology of smooth convergence on compacts.

When $\psi$ is singular at $p$, the above has no sense because $\wps$ is not continuous or even it is not defined at $0$. But there are still a very general situations where $\psi$ is singular at $p$ but we have a nice limit. This will be the case in the proof of Theorem \ref{thI}, which we shall see in section \ref{singI}.

 \subsection{Blow-up of the hypersurfaces of the flow}
 
 If $p$ is a blow-up point, there is a $x\in M$ such that $F(x,t)$ converges to $p$ when $t\to T$. Then $F^{-1}(B_R^\og({p}))$ is an open set in $M\times [0,T[$ that contains $\{x\}\times [t_0,T[$ for some $t_0\in [0,T[$. For every $t\in [t_0,T[$, let $M_t$ be the connected component of $M$ containing $x$ and contained in $F_t^{-1}(B_R^\og({p}))$, where $F_t := F(\cdot, t)$. From the choice of $R$ all the $g_t$-geodesics of $M_t$ starting from $x$ are defined for all the values of its length-arc parameter or they stop just at the boundary of $M_t$.
 
 We define now the rescaled flows $F^j(\cdot,\tau(t))$ from $M_t$ into $(\overline B_{R}^\og,\og^j)$ using the above rescaling of the metric, with $t_j>t_0$,  and the following rescaling of time
 \begin{align}\nn
 \tau = \lambda_j^2 (t-t_j) = \frac{C (t-t_j)}{T-t_j}, \qquad \tau\in\left[\frac{C (t_0-t_j)}{T-t_j},C \right[,
 \end{align}
 that is, $F^j:\bigcup_{\tau\in \left[\frac{C (t_0-t_j)}{T-t_j},C \right[} M_{t(\tau)}\times\{\tau\}
\flecha (B_R^\og, \og^j)$, with $t(\tau)= t_j+(1/C)(T-t_j) \tau = t_j+\lambda_j^{-2} \tau$ is defined by
\begin{align}\lb{reF}
 F^j(\cdot,\tau) =   F(\cdot, t_j + \lambda_j^{-2} \tau)  \text{ in } (B_R^\og, \og^j)
 \end{align}
 
 {\it For simplicity, and without loss of generality, from now on we shall take $t_0=0$.}
 
  From  $F(x,t)\underset{s\to T}{\flecha}p$, taking $t\to T$ in \eqref{dFF}, it follows that 
$d_\og(F(x,s),p) \le 2 \sqrt{C} \(\sqrt{n}+\sqrt{b}\) \sqrt{T-s}$. Then 
\begin{align}\lb{pblowup}
d_{\overline{g}^j}(F^j(x,\tau),p) \le 2 \sqrt{C} \(\sqrt{n}+\sqrt{b}\) \sqrt{\lambda_j^2(T-t_j+t_j-s)} = 2 \sqrt{C} \(\sqrt{n}+\sqrt{b}\) \sqrt{(C-\tau)}.
\end{align}  
 
 Let $N_\tau^j$ be the $\og^j$-unit normal vector of the immersion $F^j(\cdot,\tau)$ and $\{e_{\tau\ i}^j\}_{i=1}^n$ a local orthonormal tangent frame of the immersion $F^j(\cdot,\tau)$. Obviously, $N_\tau^j = N_\tau /\lambda_j$, $e_{\tau\ i}^j= e_{\tau\ i}/\lambda_j$, and the second fundamental forms $\a$ and $\a^j$ of $F$ and $F^j$ are related by:
 \begin{align}
 h_\tau(X,Y) = \og(\ona_XY, N_\tau)= \frac1{\lambda_j^2}\og^j( \ona^j_XY, \lambda_j N_\tau^j) = \frac1{\lambda_j} h_\tau^j(X,Y), \lb{al1}
 \end{align}
 which gives, for the mean curvature,
 \begin{align}
 H_\tau = \sum_i h_\tau(e_{\tau\ i},e_{\tau\ i})= \sum_i \frac1\lambda_j h_\tau^j(\lambda_j e_{\tau\ i}^j, \lambda_j e_{\tau\ i}^j) = \lambda_j H_\tau^j \lb{al2}
 \end{align}
 and, for the norm of the second fundamental form,
 \begin{align}\lb{al3}
 |h_\tau|_{\og}^2 =  \sum_{ik}|h_\tau(e_{\tau\ i},e_{\tau\ k})|^2 =\sum_{ik} \frac1{\lambda_j^2} \lambda_j^4 |h_\tau^j(e_{\tau\ i}^j,e_{\tau\ k}^j)|^2 = \lambda_j^2   |h_\tau^j|_{\og^j}^2
 \end{align}

 Moreover
 \begin{align}\lb{al4}
 \og(\ona\psi ,N)&=  \frac1{\lambda_j^2} \og^j(\lambda_j^2 \ona^j\psi , \lambda_j N^j) =  \lambda_j \og^j(\ona^j\psi , N^j)
 \end{align}

  \begin{prop}\lb{bufl} Let $F:M\times [0,T[ \flecha \oM$ be a type I \DMCF of a compact manifold $M$ closed or with boundary $\partial M$ such that $F(\partial M, t)\subset \partial G$ for a domain $G$ of $\oM$ with smooth boundary. Let $p\in G\setminus \partial G$ be a blow-up point of the flow. Let us suppose that either
  \begin{itemize}
  \item[(i)] the density $\psi$ is regular at $p$, or
  \item[(ii)] $\psi$ is singular at $p$ and the set $S$ of singular points of $\psi$ is a regular submanifold of $\oM$ satisfying that, given any curve $c:[0,1]\to \oM$ with $c(1)\in S$ and $c([0,1[)\cap S=\emptyset$, the $\lim_{t\to 1}\ona\psi/|\ona\psi|(c(t))$ lies in the normal bundle of $S$.
  \end{itemize}
Then, in case (i) the blow-up sequence of maps $F^j$ defined in \eqref{reF} subconverges smoothly on compacts to a \DMCF $\wF^\infty:M_\infty\times]-\infty,C[ \flecha\re^{n+1}$ with the Euclidean metric $g^e$ and a density $\wps^\infty$ which is constant. In case (ii), may be the limit $\wps^\infty$ is not well defined, but still the maps $F^j$ subconverge smoothly on compacts to a flow $\wF^\infty:M_\infty\times]-\infty,C[ \flecha\re^{n+1}$ with the Euclidean metric $g^e$, each embedding $\wF^\infty(\cdot,\tau)$ has a well defined  $\wps  ^\infty$-mean curvature and  $\wF^\infty(\cdot,\tau)$ follows a \DMCF motion driven by this $\wps  ^\infty$-mean curvature.
  
 In both cases the flows $F^j$ and $\wF^\infty$ are of type I, and every hypersurface $\wF^\infty(M_\infty,\tau)$ with the metric induced by the immersion $\wF^\infty(\cdot, \tau)$ is complete.
 \end{prop}
 \begin{demo} 
 From their definition and formulae \eqref{al2} and \eqref{al4}, the  $F^j(\cdot, \tau)$  satisfy the equation:
 \begin{align}\label{Fjtau}
 \parcial{F^j}{\tau}= \lambda_j^{-2}  \parcial{F}{t} = \frac{1}{\lambda_j^2} H_\psi N = \frac{1}{\lambda_j^2} (\lambda_j H^j - \lambda_j\og^j(\ona^j\psi, N^j))\  \lambda_j N^j = (H^j - \og^j(\ona^j\psi, N^j))\   N^j.
 \end{align}
Then every $F^j(\cdot,\tau)$ is a \DMCF in the ambient space with density $(\oM, \og^j, \psi)$. For every $\tau$, because $F$ is of type I, using \eqref{al3}, \eqref{al4} and Definition \eqref{defTI}, we obtain
\begin{align}\lb{tIi}
|\a_\tau^j|_{\og^j}^2 + \frac1b \og^j(\ona^j\psi,N^j)^2 = \frac1{\lambda_j^2} |\a_\tau|_{\og}^2 + \frac1{\lambda_j^2} \frac1b\og(\ona\psi,N)^2 \le \frac1{\lambda_j^2}  \frac{C}{T-t} = \frac{C}{C-\tau} .
\end{align}
Then, every flow $F^j$ is of type I, because $C$ is the supremum of the values of $\tau$.

Now, let us consider the flows $F^j(\cdot, \tau)$ for $\tau$ defined only on the closed interval $[-\lambda_{j_0}^2 t_{j_0}, C-\eps]$, $j\ge j_0$. On this interval, \eqref{tIi} gives an universal bound 
\begin{align}\lb{unib}
|\a_\tau^j|_{\og^j}^2 + \frac1b \og^j(\ona^j\psi,N^j)^2 \le \fracc{C}{C-\tau} \le \fracc{C}{C-(C-\eps)} =\fracc{C}{\eps}.
\end{align} 

\noindent Since $\og^j(\ona^j\psi,N^j)^2$ is bounded, we have only three posibilites:

 \begin{description}
\item[1] $\psi$ is smooth everywhere.
\item[2] $F^j(\cdot,\tau)$ does not touch any singular point of $\psi$ for $\tau\le C-\eps$
\item[3] There is a first $\tau_0\le C-\eps$ and a singular point $F^j(x_0,\tau_0)$ of $\psi$, such that \\ 
$\lim_{(x,\tau)\to (x_0,\tau_0)} \og^j(\ona^j\psi,N^j)(x,\tau)=0$.
\end{description}
In case 3, from the hypothesis $\lim_{z\to p} \ona^j\psi/|\ona^j\psi|$ is in the normal bundle of the singular set $S$ of $\psi$, one deduces that $F^j(M_{\tau_0},\tau_0)$ is transversal to $S$ at $F^j(x_0,\tau_0)$. Since the original hypersurface did not touch $S$, this implies that there is a $t(\tau)< t(\tau_0) < T$ where the hypersurface $F(M,t(\tau))$ is tangent to $S$, which gives $|\<\ona \psi,N\>| = \infty$ at that point, which is impossible because $T$ is the first singular time for the flow $F$. That means that case 3 is impossible. But in cases 1 and 2 the \DMCF $F^j$ in $(\oM,\og^j)$ is equivalent to the MCF $F^j\times Id$ in $\oM\times_{e^\psi}S^1$, and the norms of the second fundamental form  of the corresponding immersions satisfy $|\ha_\tau^j|^2 = |\a_\tau^j|_{\og^j}^2 + \og^j(\ona^j\psi,N^j)^2\le \frac{C(1+b)}{\eps}$. Then, the usual computations (see \cite{Hu86}) for the MCF give that all $|\hna^{jr}\ha^j_\tau|$ are bounded for every $r\in\ene$. From \eqref{hAA} and the rules for the covariant derivative in a warped product, this implies that $|\nabla^{jr}\a^j_\tau|$ and $\left|\nabla^{jr}\<\ona^j\psi,N^j\>\right|$ are bounded for every $r\in \ene$.
These bounds and \eqref{Fjtau} give also universal bounds on the derivatives of $F$ respect to $\tau$. All these bounds are taken with the metrics $\og^j$ in the ambient spaces.
 
Let us denote by  $\ona^e$ the covariant derivative in $\re^{n+1}$, and by $g_\tau^{je}$, $\wna^{je}_\tau$, $\wN_\tau^{je}$, $\wa_\tau^{je}$ and $\wH_\tau^{je}$ the  corresponding metric, unit normalal, second fundamental form and mean curvature of the immersion $\wF^j(\cdot, \tau):=\p_j^{-1}\circ F^j(\cdot, \tau)$ into $B_{\lambda_j R}^{e}$ with the Euclidean metric $g^e$. Since $(B_R^\og(p),\og^j)$ converges to $(\re^{n+1}, g^e)$ when $j\to \infty$ (that is, the metrics $\ds\{\wg^j = \p_j^*\og^j\}_{j\ge j_0}$ converge to $g^e$ on each $B_{\lambda_{j_0} R}^{g^e}$), one has that 
 \begin{align}\lb{limeuc}
 &\lim_{j\to \infty} |g_\tau^{je}- \p_j^{*}g_\tau^{j}|_{g^e} =0, \quad \lim_{j\to \infty} |\wN_\tau^{je}-\p_{j*}^{-1} N_\tau^{j}|_{g^e} =0, \quad \lim_{j\to \infty} |\wa_\tau^{je}-\p_j^{*}\a_\tau^{j}|_{g^e} =0, \nn \\
 & \lim_{j\to \infty} |\wH_\tau^{je}- H_\tau^{j}\circ\p_j|_{g^e} =0, \quad \lim_{j\to \infty} |{(\wna_\tau^{je})^k} \wa_\tau^{je}- \p_j^{*}{(\nabla_\tau^{j})^k} \a_\tau^{j}|_{g^e} =0.
 \end{align}
From \eqref{limeuc} it follows that when we consider all the above magnitudes in $\re^{n+1}$ with the Euclidean metric, they are also bounded. Moreover, by \eqref{pblowup}, the distance from $\p_j^{-1}\circ F^j\(\bigcup_{\tau\in \left[\frac{C (t_0-t_j)}{T-t_j},C \right[} M_{t(\tau)}\times\{\tau\}\)$ to $0\in \re^{n+1}$ is bounded. Then, by standard arguments (like in \cite{man} page 91 or \cite{zhu} page 87), the maps $\p_j^{-1}\circ F^j$ $C^\infty$-converge on the compacts to
 a smooth map $\wF^\infty$ defined on $M_\infty\times ]-\infty, C [$, for some manifold limit $M_\infty$, with values in $\re^{n+1}$ with the Euclidean metric (and with a density $\wps^\infty= \lim_{j\to\infty} \wps^j$ if the limit exists). 

Although may be that $\wps^\infty$ is not always well defined, there is always a limit of the mean curvatures associated to the densities, which we shall still name the mean curvature associated to the limit density. In fact, from \eqref{tIi} it follows that both $H_\tau^j$ and $\og^j(\ona^j\psi, N^j)$ are bounded, then there is a subsequence of $F^j$ such that $\lim_{j\to\infty}H^j_{\psi\tau} = \lim_{t\to\infty}(H^j_\tau-\og^j(\ona^j\psi, N^j)) =:\wH^\infty_\tau$ exists, and it is the claimed mean curvature associated to the limit density and is bounded by $\sqrt{\fracc{(n+b) C}{\eps} }$ for $\tau\le C-\eps$.

On the other hand, since we have a $C^\infty$-convergence, also respect to $\tau$,
we have
\begin{align}
\parcial{\wF^\infty}{\tau} = \lim_{j\to \infty} \parcial{\p_j^{-1}\circ F^j}{\tau} = \lim_{j\to\infty} (H^j - \og^j(\ona^j\psi, N^j))\  N^j = \wH^\infty_\tau   N^\infty .
\end{align}
  Moreover, the limit flow $\wF^\infty$ is also of type I because of \eqref{tIi} and \eqref{limeuc}.
 
 Because we have chosen $M_t$ as the connected component of $M\times \{t\}$ containing $(x,t)$ and contained in $F^{-1}(B_R^\og({p}))$, the $g^t$-geodesics of $M_t$ starting from $(x,t)$ are well defined until they touch the boundary of $M_t$ which is contained in the boundary of $B_R^\og$. But, when $j\to \infty$, this boundary goes to the infinite, then, in the limit, the geodesics starting from $(x,t)$ are well defined until the infinite. Then the corresponding limit manifold is complete.  
 \end{demo}

 \section{\DMCF of curves in surfaces producing type I singularities}\lb{singI}

\subsection{On the setting and the hypotheses of Theorem \ref{thI}}

In this section we shall describe in detail the setting for Theorem \ref{thI} and will motivate its hypotheses, then we describe some basic properties of the \DMCF in this setting, as the evolution of the barrier lines, which implies the production of singularities, and the preservation of the sign of $\kappa_\psi$ and the property of being a graph (Propositions \ref{finiteT} and \ref{Hppre} and Corollary \ref{coro9}). Then we state a serie of formulae of variation with the goal of proving (in the next section) that all the singularities which are formed are of type I.
 
We shall consider a surface $\oM$ with metric $\og = dr^2 + e^{2\p(r)} dz^2$. We consider on it a density $\psi$ which depends only on the coordinate $r$ and is singular on the geodesic $r=0$, which we shall denote also by $\Gamma$. Let us remark that the coordinate $r$ of a point gives the distance of this point to $\Gamma$.

We take $z$ as the arc-length parameter of the geodesic $r=0$ on $\oM$, then $\p$ must satisfy  $\p(0)=0$. Moreover,
 the fact that $\og$ is a metric imposes that $\p'$ has the Taylor expansion 
 \begin{align}\lb{aeph}
\p'(r) = - \overline K(0) \ r + \overline K(0)^2 \ r^3/6 + ... 
\end{align}
and the Gauss curvature $\overline K$ of $\oM$ is given by
\begin{align}\lb{oK}
\overline K = - \p''-{\p'}^2.
\end{align}
and, for the covariant derivative $\ona$ of $\oM$, we have
\begin{align}\lb{onaM}
\ona_{\partial_r} \partial_r =0, \qquad \ona_{\partial_z}\partial_z = - \p' e^{2\p} \partial_r, \qquad \ona_{\partial_r}\partial_z = \ona_{\partial_z}\partial_r = \p' \partial_z.
\end{align}

It follows from the first equation \eqref{onaM} that  the curves $z=$constant are geodesics and, from  \eqref{onaM} and \eqref{aeph}, that $\Gamma$ (curve $r=0$) is again a geodesic. Moreover, it is immediate from the expression of $\og$ that  the reflection respect to $\Gamma$ ($(r,z)\mapsto (-r,z)$ and those respect to the curves $z=c$  ($(r,c-z)\mapsto (r,c+z)$) are isometries.

From \eqref{onaM}, we obtain the following concrete expression when $M$ is a curve in a surface $\oM$ like the above one:

\begin{align}\lb{hesspsi}
\<\ona_{N} \ona\psi, N\> &= (N\psi') \<\ona r, N\> + \psi' \<\ona_N\ona r, N\> = \psi'' \<\ona r, N\>^2 + \psi' \p' \<N-\<N,\ona r\>\ona r, N\> \nn \\
& = \psi'' \<\ona r, N\>^2 + \psi' \p' \(1 -\<N,\ona r\>^2\).
\end{align} 

When $M$ is a graph $(r(z),z)$ over $\Gamma$, one has also the following useful formulae for the unit tangent vector $\bt$ and the unit normal $N$ to the curve $M$:
\begin{align}\lb{tN}
\bt = \frac{\dot r \ona r + \partial_z}{\sqrt{{\dot r}^2 + e^{2\p}}}, \qquad N= \frac{- e^{2\p} \ona r + \dot r \partial_z}{e^\p \sqrt{{\dot r}^2 + e^{2\p}}}
\end{align}
and the following expressions for its curvature $\kappa$ and for $u:=\<N,\ona r\>$
\begin{align}\lb{kappau}
\kappa = \<\ona_{\bt}\bt, N\>  = \frac{e^\p}{\sqrt{\dot{r}^2 + e^{2\p}} }\left( \frac{-\ddot{r}  + \dot{r}^2  \p'}{\dot{r}^2 + e^{2\p}} +  \p'\right), \qquad u:=\<N,\ona r\> = \frac{ - e^\p}{\sqrt{{\dot r}^2+ e^{2\p}}}.
\end{align}

When $\oM \times_{e^{\psi/m}} S^m$ is a Riemannian manifold,  for every $p\in \Gamma$, let us consider the hypersurface of $\oM \times_{e^{\psi/m}} S^m$ given by $\{\exp_p r u,\ u\in T_{p}\Gamma^\bot \subset T_p(\oM \times_{e^{\psi/m}} S^m)\}$. The Weingarten map $A^S$ of a geodesic sphere defined in that hypersurface is given by 
\begin{align*}
A^S X = -\ona_X\ona r + \<\ona_X\ona r,\frac{\partial_z}{e^\p}\>\frac{\partial_z}{e^\p} = \frac{\psi'}{m}\ X, \text{ that is, } A^S = \frac{\psi'}{m}\ Id.
\end{align*}  
From this and the properties of the mean curvature of a geodesic sphere around a point in a Riemannian manifold (see \cite{ChVa82}, Theorem 3.2), one gets
\begin{align}\lb{daspsip}
\psi'(r) = \frac{m}r - \frac{r}{3} \Ric^S(\ona  r,\ona r)(p) + \sum_{j\ge 2} A_j r^j,
\end{align}
where the $A_j$ are constants determined by the value at $p$ of polynomials in the curvature and the covariant derivatives of the curvature of $\{\exp_p r u,\ u\in T_{p}\Gamma^\bot\}$ at $p$. The condition \eqref{daspsip} implies
\begin{align}\lb{limps}
&\lim_{r\to 0}\fracc{\psi^{(n)}(r)}{m/r^n}= (-1)^{n-1} (n-1)! \ \text{ for every } n\\
\text{ and }\ &\limsup_{r\to 0}\(\frac{\psi'''}{\psi'} - \frac2{b^2} {\psi'}^2\) \quad \text{ is bounded from above}.\lb{limpsb}
\end{align}

Moreover the sectional curvatures of $\oM\times_{e^{\psi/m}}S^m$ corresponding to the planes $rz$ and those generated by $\partial_r$ and a vector $\partial_i$ tangent to $S^m$, or by $\partial_z$ and $\partial_i$ are, respectively,
\begin{align}\lb{mscur}
S_{rz}&=\widehat{R}_{rzrz}=-(\p''+{\p'}^2),\qquad S_{ri}=\widehat{R}_{riri} = -\frac{{\psi'}^2+m\psi''}{m^2},
\qquad
S_{zi}=\widehat{R}_{zizi}= -\dfrac{\psi'\p'}{m}.
\end{align}

The formulae \eqref{limps} and \eqref{limpsb} have motivated us to consider densities  $\psi$ over $\oM$ satisfying \eqref{limps2} and \eqref{limps3}. Observe, however, that, when $b$ is not a natural number, there is no natural $m$ for which $\oM \times_{e^{\psi/m}} S^m$ to be a smooth Riemannian manifold, because, in these cases, the sectional curvature $S_{ri}$ becomes infinite when $r\to 0$.

According to formulae \eqref{mscur}, when $\oM \times_{e^{\psi/m}} S^m$ is a manifold (that is $b=m$), $S_{ri}\ge 0$ is equivalent to  $\psi'' + \psi'^2/m \le 0$, then the hypothesis $r\le\sup\{r; \psi''+\psi'^2/b\le 0\}$ in Theorem \ref{thI} is not a rare analytic condition, it is motivated by this fact and can be understood as the positivity of the sectional curvature in some strange manifolds with  not integral dimension.

\subsection{The proof of points 1 to 3 of Theorem \ref{thI}}
In this subsection  we shall write many evolution formulae in an appropriate way to apply maximum principles. In those expressions will appear the laplacian $\Delta_\psi$ associated to a density,  and we recall here its definition 

\begin{align}\lb{Deltapsi}
\Delta_\psi f = \Delta f + \< \nabla \psi, \nabla f\>
\end{align}
and the way that the divergence theorem applies: Given an $n+1$ dimensional oriented compact Riemannian manifold $\oO$  with smooth boundary $\partial \oO$, let $N$ be the unit vector normal to $\partial\oO$ pointing outward, one has 
\begin{align}
&\int_{\oO} f \oDelta_\psi f \ dv^{n+1}_\psi=  - \int_{\oO} |\ona f|^2 \ dv^{n+1}_\psi + \int_{\partial {\oO}}  f \<\ona f, N\> dv^n_\psi. \lb{divtp2}
\end{align}

 \begin{nota}\lb{cilin} The following properties of the \DMCF \eqref{gmcf} will be used from now on.
 
\begin{fleqn}
\begin{equation}\lb{varrt}
\text{(a) }  \dfrac{\partial r}{\partial t}=\<\dfrac{\partial F}{\partial t},\ona{r}\>=
\<\kp N,\ona{r}\>=\kp u \quad \text{ where } \quad u:=\<N,\ona r\>.
\end{equation}
\end{fleqn}
  (b) If the evolving curve is a graph, $\psi'\ge 0$  and $N$ points to the singular axis,  $-\<\ona\psi, N\> = -\psi' \<\ona r, N\> = -\psi'\ u \ge 0$. \\
 (c) On a line $r=$constant, the value of $\kp$ is
 $ \p' -\<\ona \psi, N\> = \p' - \psi' \<\ona r, - \ona r\> = \p' + \psi'$ \\
 which is positive (in fact it is $+ \infty$ if $\lim_{r\to 0} \psi'(r)= +\infty$) at $r=0$ and remains to be positive for $r$ in the interval $[0,\z(\p' + \psi')[$.\\
 (d) If we apply the evolution formula \eqref{varrt} to the \DMCF of the line $r=r_0 < \z(\p'+\psi')$, we obtain that it evolves giving lines  $r=r(t)$ satisfying the equiation
 \begin{align}\lb{drt42}
 \parcial{r}{t} = -\p' -\psi' < - \min\{(\p'+\psi')(r), \ 0<r\le r_0\} =: - \mu , \text{ with } r(0)= r_0, \ \mu>0,
 \end{align}
 whose solution satisfies 
 $r(t)< r_0 - \mu t$, then in a finite time $T<r_0/\mu$, $r(T)=0$ and the curve $r=0$ has $\psi$-curvature $\kappa_\psi=\infty$.\\
 (e) For the \DMCF \eqref{gmcf} one has the following variational formulae (cf. \cite{mv1}), 
 
 \bec\label{evol.N}
 \frac{\ona N_t}{\partial t} = - \nabla \kappa_\psi,
 \eec
\bec\label{evol.dV}
\parcial{}{t} dv_{g_t} = - \kappa_\psi\ \kappa dv_{g_t} \text{ and } \parcial{}{t} dv^n_{\psi} = - \kappa_\psi^2 dv^n_{\psi}.
\eec

\begin{align} \lb{evolFHp}
\parcial{\kappa_\psi}{t} &= \Delta_\psi \kappa_\psi + \kappa_\psi \Big(|A|^2+ (\oRic)_{NN} -\<\ona_N\ona \psi, N\>\Big)\nn \\
&=\Delta_\psi \kappa_\psi + \kappa_\psi \Big(\kappa^2+ \overline K -\(\psi'' \<\ona r, N\>^2 + \psi' \p' \(1 -\<N,\ona r\>^2\)\)\Big).
\end{align}
 \end{nota}

 As a consequence of Remark \ref{cilin} (c) and (d) and the avoidance principle, one has
 
 \begin{prop}\lb{finiteT} Any curve inside the domain bounded by the lines $r=r_0 < \z(\p'+\psi')$ and $r=0$ which moves by the \DMCF remains contained on this domain along all the motion, and the maximal time of existence of the motion is finite $T<r_0/\mu$, where $\mu$ is defined in \eqref{drt42} . 
\end{prop}

In most of the estimates that we obtain below we apply the maximum principle. This requires that the maximum or minimum are given in an interior point. Then, in the setting where the curve that moves has boundary, we should need to study the boundary separated. But since we are working with curves contained between the lines $z=b_1$, $z=b_2$, and in our ambient surfaces the maps  $(r,\ c-z) \mapsto  (r,\ c + z)$ are isometries, if, for instance, the maximum is at $z=b_2$, we can consider the symmetry respect to $z=b_2$ given by $(r,z)=(r, \ b_2+(z-b_2)) \mapsto  (r, \ b_2-(z-b_2))=(r,\ 2\ b_2 - z)$ that doubles the curve and the point $z=b_2$ becomes an interior point, to which all the arguments apply.
 
 \begin{prop}\lb{Hppre}  Let $F(\cdot,t)$ be a compact curve evolving under \DMCF in a surface with density $(\oM,\og, \psi)$ satisfying the hypotheses of Theorem \ref{thI}. If $\kp\ge 0$ but not identically $0$ on $M_0$, then $\kp>0$ for $t>0$.
 \end{prop}
 \begin{demo} From the Proposition \ref{finiteT} and the hypothesis on the bound of $r\le r_0 <\min\{z(\p'+\psi'),  \sup\{r; (\psi'' + \psi'^2/b )|_{[0,r]}\le 0\}\}$ at time $0$, one has that $r(F(\cdot,t)) \le r_0 $ for every $t\in[0,T[$. This inequality, $\overline K\ge 0$, \eqref{limps2}, \eqref{aeph} and \eqref{oK} imply that 
 \begin{align}\lb{basineq}
 \p'\le 0, \qquad \p'' \le 0, \qquad \psi'>0, \qquad \psi''<0
 \end{align}
 and, since $\psi'$ and $\psi''$ are continuous  on the interval $]0,r_0]$ with $\lim_{r\to 0}\psi' = \infty$ and $\lim_{r\to 0}\psi'' = -\infty$, there are positive real numbers $\eps, \ \delta$ such that 
 \begin{align}\lb{basineq2}
  \qquad \psi'\ge\eps>0, \qquad \psi''\le -\delta <0, \qquad \psi'\p' \le 0
 \end{align} 
 Plugging these inequalities in \eqref{evolFHp}, we obtain that this equation has the form
\begin{align}
 \parcial{\kp}{t}  = \Delta_\psi \kp + \eta \kp \quad \text{ with } \eta \ge \delta u^2 \ge 0.
 \end{align}
  Then, by the strong  maximum principle (for instance, cf. \cite{ccgg}, page 181), we get $\kp(t)>0$ for $t>0$.
 \end{demo}

\begin{coro}\lb{coro9}
 Let $F(\cdot,t)$ be a compact curve evolving under \DMCF in a surface with density $(\oM,\og, \psi)$ satisfying the hypotheses of Theorem \ref{thI}. If $\kp\ge 0$ (but not identically $0$) and  $F(\cdot,0)$ is a graph over the geodesic $r=0$, then  $F(\cdot,t)$ is a graph for every $t\in[0,T[$.
 \end{coro} 
\noindent \begin{minipage}[b]{0.7\linewidth}
\hspace{4mm}\begin{demo}
  Since $u$ is $C^\infty$ on $M\times]0,T[$, $\sup_{x\in M} u$ is Lipschitz (then continuous) on $]0,T[$. If $\sup_{x\in M} u$ changes of sign, there must be a first $t_0$ where $(\sup_{x\in M} u)(t_0)=0$. By the boundary conditions, this supremum is attained at some interior point $x_0\in M$. At this point $u=\<\ona r, N\>=0$, then, by Proposition \ref{Hppre}, the curvature $\kappa$ of $M$ satisfies $\kappa(x_0) = \kp + \psi' u >0$. But this implies that, in a neighborhood of $x_0$, the curve $M_{t_0}$ is on one side of its tangent line (which has the equation $z=constant$) and touches this tangent line only at $x_0$. Then, in this neighborhood of $x_0$, the sign of $u$ changes at $x_0$, in contradiction with the fact that the maximum of $u$ on $M_{t_0}$ is on $x_0$ and $u(x_0)=0$.  
 \end{demo}
 \end{minipage}
 \begin{minipage}[b]{0.30\linewidth}
 \includegraphics[scale=0.3]{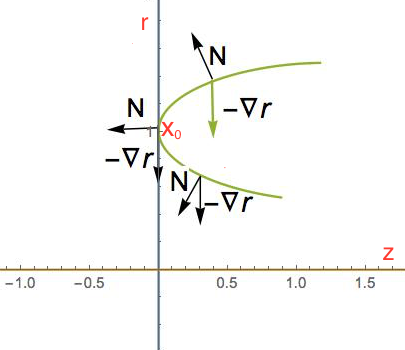}\\
 \  \\
 \  
 \end{minipage}
  
 Now we give a serie of technical lemmas (mainly variation formulae) with the aim of proving that, in this setting, the \DMCF produces Type I singularities. The strategy for doing so is similar to that used in \cite{hu90} and \cite{aag} with some variations obliged and shortcuts   possible by the circumstances of our setting.
 
 The idea is to prove that the part $k_2:= - \psi' u = -\<\nabla \psi,N\>$ (positive because $\psi'>0$ and the curve is a graph, that is, $u<0$) of $\kp=\kappa+k_2$ dominates (up to a constant) the part given by the standard curvature $\kappa$ of the curve. Then, since by hypotheses \eqref{limps2} $\psi'$ is dominated (again up to a constant) by $1/r$, the curvatures $\kappa$ and $k_2$ will be dominated by $1/r$, and it is proved at the end that satisfies the type I condition. The key points for that are good estimates of the quotients $\kappa/k_2$ and $k_2/\kappa_\psi$, which are obtained through their respective formulae of variation, whose computation requires many other formulae of variation, which we start now.
 
\begin{prop}\lb{varrD} The evolution of $r(F(\cdot,t))$ when $F(\cdot,t)$ evolves in a surface $\oM$ under the \DMCF \eqref{gmcf} is
\begin{align}
\dfrac{\partial r}{\partial t}=\Delta_\psi r
+\p'\vert\nabla r\vert^2-\p'-\psi'.
\end{align}
\end{prop}
\begin{proof}
On the curve the laplacian of $r$ is just $\bt\bt r$. Then the $\psi$-laplacian is
\begin{align*}
\Delta_\psi r
&=\mathbf{t}^{2}(r)+\<\psi' \ona r, \mathbf{t}\> \<\ona r, \mathbf{t}\>
=\mathbf{t}^{2}(r)+ \psi' \< \ona r, \mathbf{t}\>^2 
= \mathbf{t}^{2}(r)+ \psi' (1-u^2)
\end{align*}
\begin{align*}
\bt^2(r) &= \bt\<\bt, \ona r\> = \<\ona_{\bt}\bt, \ona r\> + \<\bt, \ona_\bt \ona r\> = \kappa \<\ona r,N\> + \<\bt, \p'(\bt -\<\bt, \ona r\> \ona r )\> \\
&= \kappa \<\ona r,N\> +  \p' - \p' \<\ona r,\bt\>^2 = (\kp+\<\psi' \ona r, N\>)\<\ona r, N\> +  \p' \<\ona r,N\>^2\\ 
& = \kp \<\ona r, N\> + \(\psi'+\p'\) \<\ona r,N\>^2
\end{align*}
Then
\begin{align*}
\Delta_\psi r
& = \kp \<\ona r, N\> + \psi'+\p' \<\ona r,N\>^2 \end{align*}
Plugging this expression into \eqref{varrt},
\begin{align*}
\parcial{r}{t} &= \Delta_\psi r - \psi' - \p' \<\ona r,N\>^2 
\end{align*}
what coincides with the formula that we wanted to prove.
\end{proof}

\begin{prop} The evolution of $u(F(\cdot,t))$ when $F(\cdot,t)$ evolves in a surface $\oM$ under the \DMCF \eqref{gmcf} is
\begin{align}\lb{dutlap}
\parcial{u}{t}= \Delta_\psi u -  \(\psi' \p'+ \p'^2 - \p''-\psi''\) u  (1-u^2) + \(\kappa + \p' u\)^2 u
\end{align}
\end{prop}
\begin{demo} Using \eqref{evol.N} and \eqref{hesspsi} for $\psi=r$, we obtain 
\begin{align}\parcial{u}{t} &= \kp \<\ona_N\ona r, N\> - \<\ona r, \nabla\kp\> = \kp \p'(1-u^2) - \<\ona r, \bt\>\bt\kappa  + \<\ona r, \nabla (\psi' u)\>\nn \\
&=  \kp \p'(1-u^2) - \<\ona r, \bt\>\bt\kappa  + \psi' \<\ona r, \nabla u\> + u \psi'' \< \ona r, \nabla r\> \lb{dut}
\end{align}
Computing now the Laplacian,
\begin{align*}
\Delta_\psi u &= \bt\bt u + \psi' \<\ona r, \nabla u\> \\
\bt u &= \<\ona_\bt N, \ona r\> + \<N, \ona_\bt \ona r\> = -\kappa \<\bt,\ona r\> + \<N, \p' (\bt - \<\bt, \ona r\> \ona r)\> = -\(\kappa + \p' u\) \<\bt, \ona r\>\\
\bt \bt u &= - \((\bt \kappa) + \p'' \<\bt,\ona r\> u + \p' \(-\(\kappa + \p' u\) \<\bt, \ona r\>\) \) \<\bt, \ona r\> \\
& \qquad - \(\kappa + \p' u\) \(\kappa\<N, \ona r\> + \<\bt, \p' (\bt - \<\bt, \ona r\> \ona r)\>\)\\
&= -(\bt \kappa) \<\bt, \ona r\> + \p' \kappa (1-u^2) +\(\p'^2 - \p''\) u (1-u^2) - \kappa^2 u - 2 \p' \kappa u^2  - \p'^2 u^3
\end{align*}
\begin{align}
\Delta_\psi u &= -(\bt \kappa) \<\bt, \ona r\> + \p' \kappa (1-u^2) +\(\p'^2 -\p''\) u (1-u^2) - \kappa^2 u - 2 \p' \kappa u^2  - \p'^2 u^3
\nn \\
&\qquad + \psi'\<\nabla r, \nabla u\>
\end{align}
By substitution of this expression in \eqref{dut}, we obtain:
\begin{align*}
\parcial{u}{t} &= \kp \p'(1-u^2) +\Delta_\psi u - \p' \kappa (1-u^2) - \(\p'^2 - \p''\) u (1-u^2) + \kappa^2 u + 2 \p' \kappa u^2  + \p'^2 u^3 + u \psi'' (1-u^2)\\
&= \Delta_\psi u - \(\psi' \p'+ \p'^2 - \p''-\psi''\) u(1-u^2) + \kappa^2 u+ 2 \p'\kappa u^2 + \p'^2 u^3.
\end{align*}
\end{demo}

\begin{prop} The evolution of $(\psi' u)(F(\cdot,t))$ when $F(\cdot,t)$ evolves in a surface $\oM$ under the \DMCF \eqref{gmcf} is
\begin{align}\parcial{(\psi' u)}{t} &= \Delta_\psi (\psi' u) +  \frac{-\psi''\psi' - \psi''\p'}{\psi'^3} (\psi' u)^3 + (\kappa+\p' u)^2 \psi' u\nn \\
&\quad+ \(-\psi''' u  -  \(\psi' \p'+ \p'^2 - \p''\) u \psi' + 2 \psi''\(\kappa + \p' u\) \) (1-u^2) \lb{vark2}
\end{align}
\end{prop}
\begin{demo} We just compute $\parcial{}{t} - \Delta_\psi$ acting on $\psi' u$.
\begin{align}\label{dpsiu}
\parcial{(\psi' u)}{t} - \Delta_\psi(\psi' u) = u \(\parcial{\psi'}{t} - \Delta_\psi \psi'\) + \psi' \(\parcial{u}{t} - \Delta_\psi u\) - 2 \<\nabla \psi', \nabla u\>.
\end{align}

\begin{align}
\parcial{\psi'}{t}  = \psi'' \parcial{r}{t}= \psi'' \kp u = \psi'' \(\kappa u - \psi' u^2\).
\end{align}

\begin{align}
\Delta \psi' &= \bt\bt \psi' =\bt\(\psi''\<\bt, \ona r\>\) = \psi''' (1-u^2) + \psi''\(\kappa \<N,\ona r\> + \p' \<\bt ,(\bt-\<\bt, \ona r\> \ona r)\>\) \nn \\
&= \psi'''(1-u^2) + \psi'' \kappa u + \psi'' \p' u^2.
\end{align}
\begin{align}\lb{last}
u \(\parcial{\psi'}{t} - \Delta_\psi \psi'\) &= u \(\psi'' \(\kappa u - \psi' u^2\) - \psi'''(1-u^2) - \psi'' \kappa u - \psi'' \p' u^2 - \psi'\psi''(1-u^2)\) \nn\\
&= - u \(\psi''  \psi' u^2 + (\psi'''+ \psi'\psi'')(1-u^2)  + \psi'' \p' u^2\) .
\end{align}
From \eqref{dutlap}, \eqref{dpsiu} and \eqref{last}, one obtains
\begin{align*}
\parcial{(\psi' u)}{t} - \Delta_\psi(\psi' u) &= - u \(\psi''  \psi' u^2 + (\psi'''+ \psi'\psi'')(1-u^2)  + \psi'' \p' u^2\) \\
& \quad + \psi'\(-  \(\psi' \p'+ \p'^2 - \p'' -\psi''\) u (1-u^2) + \(\kappa + \p' u\)^2 u\)\\
&\quad + 2 \psi''(\kappa + \p' u) (1-u^2)
\end{align*}
and formula \eqref{vark2} follows.
\end{demo}

 If the curve is a graph over the geodesic $r=0$ and  $\psi'\neq 0$, then $k_2:= -\psi'u >0$.

Now, we would like to prove that $k_2\leq C k_{\psi}$ for some constant $C$ independent of $t$. Let us observe that if  $\kappa\geq 0$, as $k_2>0$ and $\kappa_\psi>0$, then $k_2\leq k_{\psi}$ which gives $\dfrac{k_2}{\kappa_\psi}\leq 1$, then the difficulties to prove the bound we like arise only when $\kappa<0$. We start by computing the variation of the quotient $\dfrac{k_2}{\kappa_\psi}$.

\begin{lema}\lb{vark2kp} The evolution of $\dfrac{k_2}{\kappa_\psi}$ when $M_t$ evolves in a surface $\oM$ under the \DMCF \eqref{gmcf} is
\begin{align}
\dfrac{\partial}{\partial t}\Big(\dfrac{k_2}{\kappa_\psi}\Big)&=\Delta_{\psi}\Big(\dfrac{k_2}{\kappa_\psi}\Big)+\dfrac{2}{\kappa_\psi}\<\nabla\Big(\dfrac{k_2}{\kappa_\psi}\Big),\nabla \kappa_\psi\>
\nn\\
&\qquad+ \dfrac{k_2}{\kappa_\psi}\Bigg(\dfrac{\psi''\varphi'}{\psi'}(2-3u^2)+2u^2\varphi'^2
-\dfrac{\psi'''}{\psi'}(1-u^2)
+\varphi''(2-u^2)
\Bigg)\nn\\ 
&\qquad  + 2\(\dfrac{k_2}{\kappa_\psi}-1\)\big(\varphi'\psi'u^2+\psi''(1-u^2)\big).\lb{qk2kp}
\end{align}
\end{lema}
\begin{proof}
We just compute using the formula $\Delta_\psi\frac{k_2}{\kappa_\psi} = -\frac{k_2}{\kappa_\psi^2} \Delta_\psi \kappa_\psi + \frac{1}{\kappa_\psi} \Delta_\psi k_2 -\frac{2}{\kappa_\psi} \<\nabla\(\frac{k_2}{\kappa_\psi}\), \nabla \kappa_\psi\>$ and the equations \eqref{evolFHp} and \eqref{vark2}
\begin{align*}
\dfrac{\partial}{\partial t}\Big(\dfrac{k_2}{\kappa_\psi}\Big)&=
\Delta_{\psi}\Big(\dfrac{k_2}{\kappa_\psi}\Big)+\dfrac{2}{\kappa_\psi}\<\nabla\Big(\dfrac{k_2}{\kappa_\psi}\Big),\nabla \kappa_\psi\>
+\dfrac{1}{\kappa_\psi}\Big(\dfrac{\partial}{\partial t}k_2-\Delta_\psi k_2\Big)-\dfrac{k_2}{\kappa_\psi^2}\Big(
\dfrac{\partial}{\partial t}\kappa_\psi-\Delta_\psi k_{\psi}
\Big)
\nn\\
&=\Delta_{\psi}\Big(\dfrac{k_2}{\kappa_\psi}\Big)+\dfrac{2}{\kappa_\psi}\<\nabla\Big(\dfrac{k_2}{\kappa_\psi}\Big),\nabla \kappa_\psi\>
\nn\\
&\qquad+\dfrac{k_2}{\kappa_\psi}\Bigg(-\dfrac{\psi''(\psi'+\varphi')}{\psi'^3}k_2^2+\kappa^2+2\kappa\varphi'u+\varphi'^2u^2
\nn\\
&\quad\qquad\qquad
-\(\dfrac{\psi'''}{\psi'}+\p'(\psi'+ \p')- \p'' \)  (1-u^2)
-2 \psi''(\dfrac{\kappa}{k_2} - \dfrac{\p'}{\psi'} ) (1-u^2)\Bigg)
\nn\\
&\qquad-\dfrac{k_2}{\kappa_\psi}
\big(\kappa^2-\varphi''-\varphi'^2-\psi''u^2-\psi'\varphi'(1-u^2)\big);
\\
\dfrac{\partial}{\partial t}\Big(\dfrac{k_2}{\kappa_\psi}\Big)&-\Delta_{\psi}\Big(\dfrac{k_2}{\kappa_\psi}\Big)-\dfrac{2}{\kappa_\psi}\<\nabla\Big(\dfrac{k_2}{\kappa_\psi}\Big),\nabla \kappa_\psi\>
\nn\\
& = \dfrac{k_2}{\kappa_\psi}\Bigg(-\dfrac{\psi''\varphi'}{\psi'}u^2-2\dfrac{\kappa}{k_2}\varphi'\psi'u^2+\varphi'^2(1+u^2)
-\(\dfrac{\psi'''}{\psi'}+\p'^2- \p'' \)  (1-u^2)
-2 \psi''(\dfrac{\kappa}{k_2} - \dfrac{\p'}{\psi'} ) (1-u^2)
+\varphi''
\Bigg)   \\
& = \dfrac{k_2}{\kappa_\psi}\Bigg(\dfrac{\psi''\varphi'}{\psi'}(2-3u^2)-2\dfrac{\kappa}{k_2}\big(\varphi'\psi'u^2+\psi''(1-u^2)\big)+2u^2\varphi'^2
-\dfrac{\psi'''}{\psi'}(1-u^2)
+\varphi''(2-u^2)
\Bigg).\\
& = \dfrac{k_2}{\kappa_\psi}\Bigg(\dfrac{\psi''\varphi'}{\psi'}(2-3u^2)+2u^2\varphi'^2
-\dfrac{\psi'''}{\psi'}(1-u^2)
+\varphi''(2-u^2)
\Bigg) -2\(1-\dfrac{k_2}{\kappa_\psi}\)\big(\varphi'\psi'u^2+\psi''(1-u^2)\big).\\
\end{align*}
Where, in the last equality, we have used that $\fracc{\kappa}{\kappa_\psi} = 1-\frac{k_2}{\kappa_\psi}$.
\end{proof}

\begin{lema}\lb{bk2kp} Under the hypotheses of Theorem \ref{thI}, for any $t_0\in ]0,T[$,  the quotient $\dfrac{k_2}{\kappa_\psi}$ is uniformly bounded on $M\times [t_0,T[$ by $\max\{ 1, \Big(\max_{M_{t_0}}\dfrac{k_2}{\kappa_\psi}\Big)e^{\alpha r_0/\mu}\}$, where $r_0 = \max \{r(x), x\in M_{t_0}\}$ and $\mu = \min\{(\p'+\psi')(r), \ 0<r\le r_0\}$.
\end{lema}
\begin{demo} When the hypotheses of Theorem \ref{thI} are satisfied, we know from Proposition \ref{Hppre}  that $\kappa_\psi(t)>0$ for every $t\in ]0,T[$, and the quotient $\dfrac{k_2}{\kappa_\psi}$ is well defined for such $t$. Moreover, from Corollary \ref{coro9} and the first equality in \eqref{basineq2}, $\dfrac{k_2}{\kappa_\psi}>0 $.
From the inequalities \eqref{basineq} it follows that  $\varphi'\psi'u^2+\psi''(1-u^2)\le 0$. Moreover, it follows from \eqref{limps2} that $-\fracc{\psi'''}{\psi'}(1-u^2) \le C_2$ for some $C_2>0$ and every $r\in ]0,r_0]$ and, from \eqref{limps2} and \eqref{aeph}, that $\fracc{\psi'' \p'}{\psi'}(2- 3 u^2) \le C_1$ for some $C_1>0$ and every $r\in ]0,r_0]$ and obviously $2 u^2 \p'^2 \le C_3>0$ and $\p''(2-u^2)\le 0$ by \eqref{basineq}. Plugging these inequalities in \eqref{qk2kp} we obtain that either $\fracc{k_2}{\kappa_\psi} \le 1$ or
\begin{align*}
\dfrac{\partial}{\partial t}\Big(\dfrac{k_2}{\kappa_\psi}\Big)&\leq \Delta_{\psi}\Big(\dfrac{k_2}{\kappa_\psi}\Big)+\dfrac{2}{\kappa_\psi}\<\nabla\Big(\dfrac{k_2}{\kappa_\psi}\Big),\nabla \kappa_\psi\> + \dfrac{k_2}{\kappa_\psi}\Bigg(C_1+ C_3 + C_2 \Bigg) 
\end{align*}
By the maximum principle, $k_2/\kappa_\psi$ is bounded from above by the solution of the equation $y'(t)= \a\ y(t)$, $\a = C_1+C_2+C_3$, with the initial condition $y(t_0)=\max_{M_{t_0}}(k_2/\kappa_\psi)$, that is
\begin{align*}
\dfrac{k_2}{\kappa_\psi}(p,t)\leq \Big(\max_{M_{t_0}}\dfrac{k_2}{\kappa_\psi}\Big)e^{\alpha (t-t_0)},
\end{align*}
for every $t\in[t_0,T[$. By Remark \ref{cilin}, $T$ is finite and lower than $r_0/\mu$, then 
\begin{align}
\dfrac{k_2}{\kappa_\psi}(p,t)\leq \Big(\max_{M_{t_0}}\dfrac{k_2}{\kappa_\psi}\Big)e^{\alpha (T-t_0)} \le \Big(\max_{M_{t_0}}\dfrac{k_2}{\kappa_\psi}\Big)e^{\alpha r_0/\mu} .
\end{align}
The statement of the Lemma follows from these remarks.
\end{demo}

\begin{lema}\lb{bkk2} For every real number $b> 0$, one has
\begin{align}
\dfrac{\partial}{\partial t}\Big(\dfrac{\kappa}{k_2}\Big)&=\Delta_{\psi}\Big(\dfrac{\kappa}{k_{2}}\Big)
+\dfrac{2}{k_{2}}\<\nabla\Big(\dfrac{\kappa}{k_{2}}\Big),\nabla k_{2}\>+2\dfrac{\kappa_\psi}{k_2}\Big[
\big(\dfrac{b\ \psi''+\psi'^2}{b}(1-u^2)+\varphi'\psi'u^2\big)\dfrac{\kappa}{k_2}
+\dfrac{1}{b}\psi'^2(1-u^2)\Big(\dfrac{1}{b}-\dfrac{\kappa}{k_2}\Big)\Big]
\nn\\
&\qquad+\dfrac{\kappa_\psi}{k_{2}}\Bigg[
-\varphi''(2-u^2)
-2\varphi'^2u^2
+\Big(\dfrac{\psi'''}{\psi'}-\dfrac{2}{b^2}\psi'^2\Big) (1-u^2)
+ \dfrac{\psi''\p'}{\psi'}(-2+3u^2)
\Bigg].
\end{align}
\end{lema}
\begin{proof} Computing like in the proof of Lemma \eqref{vark2kp}, we obtain
\begin{align*}
\dfrac{\partial}{\partial t}\Big(\dfrac{\kappa}{k_{2}}\Big)&=
\Delta_{\psi}\Big(\dfrac{\kappa}{k_{2}}\Big)
+\dfrac{2}{k_{2}}\<\nabla\Big(\dfrac{\kappa}{k_{2}}\Big),\nabla k_{2}\>
+2\big(\psi''(1-u^2)+\varphi'\psi'u^2\big)\dfrac{\kappa_\psi}{k_2}\dfrac{\kappa}{k_2}
+\dfrac{2}{b^2}\psi'^2\dfrac{\kappa_\psi}{k_2}(1-u^2)
\nn\\
&\qquad+\dfrac{\kappa_\psi}{k_{2}}\Bigg[
-\varphi''(2-u^2)
-2\varphi'^2u^2
+\Big(\dfrac{\psi'''}{\psi'}-\dfrac{2}{b^2}\psi'^2\Big) (1-u^2)
+ \dfrac{\psi''\p'}{\psi'}(-2+3u^2)
\Bigg]\nn\\
&=\Delta_{\psi}\Big(\dfrac{\kappa}{k_{2}}\Big)
+\dfrac{2}{k_{2}}\<\nabla\Big(\dfrac{\kappa}{k_{2}}\Big),\nabla k_{2}\>+2\dfrac{\kappa_\psi}{k_2}\Big[
\big(\psi''(1-u^2)+\varphi'\psi'u^2\big)\dfrac{\kappa}{k_2}
+\dfrac{1}{b^2}\psi'^2(1-u^2)\Big]
\nn\\
&\qquad+\dfrac{\kappa_\psi}{k_{2}}\Bigg[
-\varphi''(2-u^2)
-2\varphi'^2u^2
+\Big(\dfrac{\psi'''}{\psi'}-\dfrac{2}{b^2}\psi'^2\Big) (1-u^2)
+ \dfrac{\psi''\p'}{\psi'}(-2+3u^2)
\Bigg]  \nn\\
&=\Delta_{\psi}\Big(\dfrac{\kappa}{k_{2}}\Big)
+\dfrac{2}{k_{2}}\<\nabla\Big(\dfrac{\kappa}{k_{2}}\Big),\nabla k_{2}\>+2\dfrac{\kappa_\psi}{k_2}\Big[
\big(\dfrac{b\ \psi''+\psi'^2}{b}(1-u^2)+\varphi'\psi'u^2\big)\dfrac{\kappa}{k_2}
+\dfrac{1}{b}\psi'^2(1-u^2)\Big(\dfrac{1}{b}-\dfrac{\kappa}{k_2}\Big)\Big]
\nn\\
&\qquad+\dfrac{\kappa_\psi}{k_{2}}\Bigg[
-\varphi''(2-u^2)
-2\varphi'^2u^2
+\Big(\dfrac{\psi'''}{\psi'}-\dfrac{2}{b^2}\psi'^2\Big) (1-u^2)
+ \dfrac{\psi''\p'}{\psi'}(-2+3u^2)
\Bigg].
\end{align*}
\end{proof}

\begin{lema}\lb{k/k2} Under the hypotheses of Theorem \ref{thI},   the quotient  $|\kappa/k_2|$ is uniformly bounded on $M_t$ for  $t\in [0, T[$
\end{lema}
\begin{demo} If $\kappa/k_2 \le 1/b$, $\kappa/k_2$ is bounded from above. If $\kappa/k_2 > 1/b$, then in the formula of Lemma \ref{bkk2}, the addend that contains the term $1/b - \kappa/k_2$ becomes negative, $\kappa_\psi/k_2 >0$, and, thanks again to \eqref{limps2}, \eqref{basineq} and \eqref{basineq2}, the coefficient of $\kappa/k_2 $ is non positive, and the other addend that multiplies  $k_ \psi/k_2 = 1 + \kappa/k_2$ is bounded by some constant $\delta$. Then we can write 
\begin{align}
\parcial{}{t}\(\frac{\kappa}{k_2}\) &\le  \Delta_{\psi}\Big(\dfrac{\kappa}{k_{2}}\Big)
+\dfrac{2}{k_{2}}\<\nabla\Big(\dfrac{\kappa}{k_{2}}\Big),\nabla k_{2}\>+ \delta \dfrac{\kappa}{k_2} + \delta
\nn
\end{align}
and the maximum principle gives $\fracc{\kappa}{k_2} \le (1+\max_{M_0}\fracc{\kappa}{k_2}) e^{\delta T} - 1 \le (1+\max_{M_0}\fracc{\kappa}{k_2}) e^{(\delta r_0)/\mu} - 1$. Then $\kappa/k_2$ is bounded from above.

If $\kappa<0$, as by Proposition \ref{Hppre} $k_\psi >0$, we have $|\kappa| < k_2$ and  $|\kappa/k_2|\le 1$.
\end{demo}

\begin{coro}\lb{singz}
Singularities of the flow occurs when and only when the evolving curve touches the axis $\Gamma$.
\end{coro}

\begin{teor}\lb{thp3}
Under the hypotheses of Theorem \ref{thI}, the \DMCF develops, in the first singular points, singularities of type I.
\end{teor}
\begin{demo}
Instead of \eqref{gmcf}, we can use the equivalent flow 
\begin{align}\lb{egmcf}
\parcial{F}{t} = \frac{\kappa_\psi}{\<N,\ona r\>} \ona r 
\end{align}
Which has sense when the evolving curve is a graph over $\Gamma$, because this implies that $\<N,\ona r\> <0$ never vanishes, and it is equivalent to \eqref{gmcf} because $\ds \<\parcial{F}{t}, N\> = \kappa_\psi$. Under this flow, the variation of $r$ is given by
\begin{align}\lb{drdte}
\parcial{r}{t}= \<\ona r, \parcial{F}{t}\> = \< \ona r, \frac{\kappa_\psi}{\<N,\ona r\>} \ona r\> = \frac{\kappa_\psi}{u} = \frac{\kappa}{u} - \psi'.
\end{align}
From Lemmas \ref{k/k2} and \ref{bk2kp} $|\kappa| \le C_1 (-\psi'\ u) \le C_2 (-\psi'\ u + \kappa)$, then (remeber $u<0$)
$\fracc{\kappa}{u} \le \left|\frac{\kappa}{u}\right| \le C_2 (\psi'\  - \fracc{\kappa}{u})$, that is $(1+C_2) \fracc{\kappa}{u} \le C_2\ \psi'$. From this and \eqref{drdte},
\begin{align}\lb{drdte2}
\parcial{r}{t}\le \(\frac{C_2}{1+C_2} - 1\)\psi' = -\frac{1}{1+C_2} \psi' \le - C_3 \frac{b}{r},
\end{align}
where we have used \eqref{limps2} for the last inequality
then 
\begin{align}
\parcial{r^2}{t}\le  -2\ b\ C_3 =: - C_4, \nn
\end{align}
and, for any $0<t<t_1<T$, one has
\begin{align}
r^2(t_1)- r^2(t) \le   - C_4 (t_1-t), \quad \frac1{r^2(t)} \le \frac1{r^2(t_1)+ C_4 (t_1-t)} \le \frac1{ C_4 (t_1-t)} \quad \text{ for every $t_1 <T$ }\nn
\end{align}
Taking limits when $t_1 \to T$, we have the inequality
\begin{align}\lb{1/r}
\frac1{r(t)} \le \frac1{\sqrt{ C_4 (T-t)}}
\end{align}
On the other hand, by Lemma \ref{k/k2}, the definition of $k_2$, and \eqref{1/r},
\begin{align}
|k_2|^2 + |\kappa|^2  \le (1+C_5) k_2^2 \le (1+C_5) {\psi'}^2 \le C_6 \frac{b^2}{r^2} \le \frac{C_6 b^2}{C_4 (T-t)}
\end{align}
which shows that the singularity is of type I.
\end{demo}

\section{Convergence of the blow-ups for some type I singularities (the proof of point 4 in Theorem \ref{thI})}\lb{ConBU}

Let us suppose that we have a \DMCF on $\re^{n+1}$ {\it with the Euclidean metric and with a density $\psi(p) = b \ln r(p)$, where $r(p)$ is the Euclidean distance from $p$ to the axis $z\equiv x_{n+1}$ of $\re^{n+1}$}. One has the following monotonicity formula analogous to the Huisken's formula for the MCF in \cite{hu90}.

\begin{prop} Let $\uu:\re^{n+1}\times[0,T[ \flecha \re$ be defined as $\uu(p,t) = (4 \pi (T-t))^{-(n+b)/2} e^{-|p|^2/(4(T-t))}$. If $F:M\times[0,T[ \flecha \re^{n+1}$ is a family of immersions of a hypersurface $M$ moving by the \DMCF that either is compact or the Euclidean $(n-1)$-volume of the boundary of the intersections of $M_t$ with the closed balls $\overline B_R$ of $\re^{n+1}$ centered at the origin are bounded by $f(t) R^q$, with $f(t)>0$ and $q$ a fixed positive number, one has
\begin{align}\lb{formon}
\deri{}{t}\int_M \uu(F(x,t),t) dv_\psi = - \int_M \(H_\psi + \frac{\<F(x,t),N(F(x,t))\>}{2(T-t)}\)^2 \uu(F(x,t),t) dv_\psi
\end{align}
When $b=m\in\ene$, we recover the standard Huisken's monotonicity formula restricted to hypersurfaces in $\re^{n+1+m}$ obtained by the rotation of a hypersurface in $\re^{n+1}$. 
\end{prop}
\begin{demo}
First a Minkowski's formula. For hypersurfaces $X:M\flecha \re^{n+1}$ the classical Minkowski formula states that  $\Delta (\frac12 |X|^2) = n + H\< N,X\>$. When we have also a density, we have $\Delta_\psi (\frac12 |X|^2) = n + H \<N, X\> + \<\nabla \psi, \nabla (\frac12 |X|^2)\> = n + H \<N, X\> + \<\ona \psi, \frac12 \ona |X|^2\> - \<\ona\psi,N\> \<\frac12 \ona|X|^2,N\> =n + \Hp \<N,X\> + \<\ona\psi,X\>$. If $\psi(p) = b\ \ln r(p)$, $\<\ona\psi,X\>= b \frac1{r(X)} \<\ona r, X\> = b$, then
\begin{align}\lb{Mink}
\Delta_\psi (\frac12 |X|^2) = n + b + \Hp \<N,X\>.
\end{align} 
Now we take the derivative, taking into account \eqref{evol.dV}
\begin{align}
\deri{}{t}\int_M \uu(F(x,t),t) dv_\psi = \int_M \(\frac1{2(T-t)} \(n+b-\frac{|F|^2}{2(T-t)}- \Hp \<F,N\>\) - \Hp^2\) \ \uu\ dv_\psi
\end{align}
and substituting \eqref{Mink} in the above expression
\begin{align} \lb{derm}
\deri{}{t}\int_M \uu(F(x,t),t) dv_\psi = \int_M \(\frac1{2(T-t)} \Delta_\psi \(\frac12 |F|^2\)-\frac{|F|^2}{4(T-t)^2}-2 \frac{\Hp \<F,N\>}{2 (T-t)} - \Hp^2\) \ \uu\ dv_\psi
\end{align}
\begin{align}\lb{F2}
\int_M \frac{|F|^2}{4(T-t)^2} dv_\psi = \int_M \frac{|F^\top|^2}{4(T-t)^2} dv_\psi + \int_M \frac{\<N,F\>^2}{4(T-t)^2} dv_\psi
\end{align}
If $M$ is compact, we can apply divergence theorem \eqref{divtp2} for the $\psi$-laplacian with the $\psi$-volume to  
\begin{align}
\int_M \(\frac1{2(T-t)} \Delta_\psi \(\frac12 |F|^2\)\) \ \uu\ dv_\psi &= - \int_M \frac1{2(T-t)} \<\nabla \(\frac12 |F|^2\), \nabla \uu\> dv_\psi \nn\\
&=  \int_M \frac1{4(T-t)^2} \<\nabla \(\frac12 |F|^2\), \nabla \(\frac12 |F|^2\)\>  \uu dv_\psi =  \int_M \frac1{4(T-t)^2} |F^\top|^2  \uu dv_\psi \lb{psiSto}
\end{align}
By substitution of \eqref{F2} and \eqref{psiSto} into \eqref{derm}
\begin{align}
\deri{}{t}\int_M \uu(F(x,t),t) dv_\psi = \int_M\(-\frac{\<N,F\>^2}{4(T-t)^2}-2 \frac{\Hp \<F,N\>}{2 (T-t)} - \Hp^2\) \ \uu\ dv_\psi
\end{align}
which gives \eqref{formon}.

If $M$ is not compact, we can compute the integral along $M$ as the limit of the integrals along its intersections $M_{tR}:=M_t \cap\overline B_R$ with the closed balls $\overline B_R$ of $\re^{n+1}$ centered at the origin with radius $R$ when $R\to \infty$.
\begin{align}
\int_{M_{tR}} \(\frac1{2(T-t)} \Delta_\psi \(\frac12 |F|^2\)\) \ \uu\ dv_\psi &= - \int_{M_{tR}} \frac1{2(T-t)} \<\nabla \(\frac12 |F|^2\), \nabla \uu\> dv_\psi + \int_{\partial M_{tR}} \frac{1}{2(T-t)} \nu \(\frac12 |F|^2\) \uu \ dv_\psi\nn\\
& =  \int_{M_{tR}} \frac1{4(T-t)^2} |F^\top|^2  \uu dv_\psi + \int_{\partial M_{tR}} \frac{1}{2(T-t)}\nu \(\frac12 |F|^2\) \uu \ dv_\psi \lb{psiSto2} 
\end{align}
where $\nu$ is the outward unit normal vector field on $\partial M_{tR}$. Let us study  the last addend in \eqref{psiSto2}
\begin{align}
\int_{\partial M_{tR}} \nu \(\frac12 |F|^2\) \uu \ dv_\psi &= (4\pi(T-t))^{-(n+b)/2} e^{-R^2/(4(T-t))} \int_{\partial M_{tR}} \<\nu,F\> r(F)^b dv \nn \\
&\le (4\pi(T-t))^{-(n+b)/2} e^{-R^2/(4(T-t))} R^{b+1} \int_{\partial M_{tR}} dv \nn \\
&\le (4\pi(T-t))^{-(n+b)/2} e^{-R^2/(4(T-t))} f(t)\ R^{q+b+1} \to 0 \text{ when } R\to \infty.
\end{align}
and we continue the proof of the formula as in the compact case.
\end{demo}
We want to apply the above formula to a \DMCF of a complete graph $(r(z),z)$ over the axis $z$ in $\re^2$. In this case $M$ is not compact and it is easier to check  the condition of the above theorem if we take closed squares $C_R$ centered at $0$ of side $2R$ instead of balls $\overline B_R$. To check this condition we shall need the Sturmian Theorem that we shall write below.

In order to prove the announced Sturmian Theorem, we shall work with the flow \eqref{egmcf} equivalent to the \DMCF and used in the proof of Theorem \ref{thp3}. Under this flow the variation of $r$ is given by \eqref{drdte}. To obtain the variation of $\eta:= \dot r\equiv\parcial{r}{z}$ we take the derivative of  \eqref{drdte} respect to $z$ and obtain 
\begin{align}
\parcial{\dot r}{t}=   \parcial{}{z}\frac{\kappa_\psi}{u} = \parcial{}{z}\(\frac{\kappa}{u} - \psi'\) 
=\frac{u \dot\kappa - \kappa \dot u}{u^2} - \psi'' \dot r.
\end{align}

Plugging  \eqref{kappau} in this and doing the corresponding derivatives, we obtain, for $\eta := \dot r$
\begin{align}\lb{ddoteta}
\dfrac{\partial\eta}{\partial t}&=
\dfrac{\ddot{\eta}}{\eta^2+e^{2\p}}+T(\eta)\eta
\end{align}
where
\begin{align*}
T(\eta)&:=
-\dfrac{2\dot{\eta}}{(\eta^2+e^{2\p})^2}(\dot{\eta}+e^{2\p}\varphi')
-2\Big(\dfrac{\dot{\eta}}{\eta^2+e^{2\p}}-\dfrac{\eta^2}{(\eta^2+e^{2\p})^2}(\dot{\eta}+e^{2\p}\varphi')\Big)\varphi'
\nn\\
&\qquad-\dfrac{\eta^2}{\eta^2+e^{2\p}}\varphi''
-\varphi''-\psi'',
\end{align*}
From Corollary \ref{singz} we know that for every $t_0\in[0,T[$, $T(\eta)$ and also the coefficient of $\ddot\eta$ in \eqref{ddoteta} are bounded on $J\times[0,t_0]$, where $J$ is the domain where $z$ lives. Then we can apply the Sturmian Theorem of Angenent (cf. \cite{an88}) to obtain

\begin{lema}\label{discrt} Let $M_t$, with $t\in[0,T[$, be a maximal  solution of  \eqref{egmcf}  with initial condition $M_0$ and satisfying the hypotheses of Theorem \ref{thI}. For each $t\in[0,T[$, the set $Z_t=\{z\in S^1 \text{ or } z\in[a_1,a_2]; \dot{r}_t =0\}$ is finite, the function $t \mapsto N(t):=\sharp(Z_t)$ is non increasing and, at the points $(z_0,t_0)$ satisfying $0  = \dot{r}_{t_0}(z_0) =\ddot{r}_{t_0}(z_0)$ there is a neighborhood where the number of zeroes decreases.
\end{lema}

In the following Lemma, which is an adaptation of formula (57) in \cite{mv1} before doing the integration, we shall use  the following multi-index notation. Capitals will denote multi-indices. For us, all the entries $j_k$ of a multi-index  $J=(j_1, ...., j_q)$ will be ordered $j_1 \ge j_2 \ge \cdots \ge j_q >0$. For such a multi-index, we shall denote  $|J|:= j_1 + ...+ j_q$, $d(J):=q$, $o(J)=j_1$, $\partial_s^{J}x := \partial_s^{j_1}x \dots \partial_s^{j_q}x$, $\ona^{J}x := \ona^{j_1}x\otimes \dots \otimes\ona^{j_q}x$.

\begin{lema}[\cite{mv1}] One has the following evolution formula under \DMCF \eqref{gmcf} in $\re^2$,
\begin{align}
\parcial{}{t}(\partial_{s}^n \kappa_{\psi})^{2} 
&= \Delta_\psi \(\partial_{s}^{n} \kappa_{\psi}\)^2
- 2  \(\partial_{s}^{n+1} \kappa_{\psi}\)^2  \nn \\
& \qquad + 2  (a_{n0}  + a_{n1} \kappa_\psi +a_{n2} \kappa_\psi^2) (\partial_s^n \kappa_\psi)^2   + \sum a_{iJ}   \(\kappa_\psi^i \partial_s^J\kappa_\psi\) (\partial_s^n \kappa_\psi) ,
\label{E_igualdadprincipal}
\end{align}
where $i+|J|\ge 1$, $i\le n+1$, $o(J) \le n-1$, $|J| \le n$, the coefficients \lq\lq$a_{n j}$'' are polynomials of degree $2-j$ in the variables $\ona^m \psi$ (where $\ona^m\psi$ has degree $m$), and the coefficients $a_{iJ}$ are polynomials of degree $n+3-i-|J|-d(J)$ in the variables $\ona^m \psi$  acting on $\partial_s$ and/or $N$, and  some of them can be zero.
\end{lema}

Now, we are ready to prove points 4 in Theorem \ref{thI} and Corollary \ref{corThI}
\begin{teor} Under the hypotheses of Theorem \ref{thI}, at the first singular time, at each singular point, a blow-up centered at this point gives a new type I limit flow in $\re^2$ with its Euclidean metric and density $\wps^\infty = \ln r^b$ which is a graph over $r=0$ for every time and, after doing a new blow-up,  converges to a $\ln r^b$-shrinker in $\re^2$, which is the line $r=${\rm constant} in case $b=m\in \ene$.
\end{teor}
\begin{demo}
From Proposition \ref{bufl}, a blow-up centered at this point gives a limit flow in $\re^2$ with its Euclidean metric and density $\lim_{j\to\infty}\psi\circ \p_j$ is this limit exist. To know that, in fact, the limit exists and what it is, we use the property $\wps^j(v) = \wps^0(v/\lambda_j)$ stated in \eqref{wpsj}. This formula implies that  $\wps^j{}'(r) = 1/\lambda_j \wps^0{}'(r/\lambda_j)$ then, by \eqref{limps2},
$$ \lim_{j\to\infty}\frac{\wps^j{}'(r)}{b/r} =  \lim_{j\to\infty}\frac{\wps^0{}'(r/\lambda_j)}{b/(r/\lambda_j)} =1$$
that is, in  the $C^0$ convergence on compacts, there is a limit function $\wps^\infty{}'= b/r$ of $\wps^j{}'$ which allows us to define the density $\wps^\infty$ on $\re^2$ by $\wps^\infty=b \ln r$ and the limit flow  $\wF^\infty(\cdot,\tau)$ satisfies the equation \eqref{gmcf} with the mean curvature associated to this density.

Since $\wF^j(\cdot,\tau)$ are graphs  that converge $C^\infty$ on the compacts to $\wF^\infty(\cdot,\tau)$, then the equivalent flows $(\wwr^j(z,\tau),z)$ have the derivatives of $\wwr^j(z,\tau)$ respect to $z$ are bounded on every compact by the bounds of the derivatives of $\wF^j(z,\tau)$ respect to $z$, then the $\wwr^j(z,\tau)$ converge  $C^\infty$ on the compacts to a function $\wwr^\infty(z,\tau)$ and $\wF^\infty(\cdot,\tau)$ is a graph for every $\tau$.

Now, we apply to the flow $\wF^\infty(\cdot, \tau)$ in the Euclidean space with density $(\re^2,g^e,\wps^\infty=\ln r^b)$ the standard blow-up 
\begin{align}\lb{sbw1}
&\lambda(\tau)^2 = \frac{1}{2(C-\tau)}, \quad \wtau(\tau) =  \ln \lambda(\tau),
 \quad
 \wF(\cdot,\wtau) =  e^\wtau \wF^\infty(\cdot,  \tau(\wtau))   \\ \text{ and } 
 &\wps (v) = \wps^\infty(v/\lambda) = \ln{r(v/\lambda)}^b = \ln \frac{r(v)^b}{\lambda^b} , \lb{sbw1a} 
  \end{align}
which gives
 \begin{align}
 &\wka(\wF(\cdot,\wtau)):=\wH(\wF(\cdot,\wtau)) = \frac1\lambda \wH^\infty(\wF^\infty(\cdot,\tau(\wtau)))=:\frac1\lambda \wka^\infty(\wF^\infty(\cdot,\tau(\wtau))), \lb{sbw2a} \\
 & g^e(\ona\wps, N)(\wF(\cdot,\wtau)) = \frac1\lambda g^e(\ona\wps^\infty, N)(\wF^\infty(\cdot,\tau(\wtau))), \lb{sbw2aa} \\
 &|\wka|^2=|\wa_\wF|^2 =  \frac1{\lambda^2} |\wa^\infty|^2 = \frac1{\lambda^2} |\wka^\infty|^2 \lb{sbw2b} 
 \end{align}
 and, taking into account the estimate \eqref{tIi}, 
 \begin{align}
 &\frac1b g^e(\ona\wps, N)^2 +|\wka|^2 = \frac1{\lambda^2} g^e(\ona\wps^\infty, N)^2 + \frac1{\lambda^2} |\wka^\infty|^2 \le \frac1{\lambda^2} \frac{C}{C-\tau} =2\ C  \lb{sbw3}\\
 &\wka_\psi^2 = \frac1{\lambda^2} (\wka_\psi^\infty)^2 =\frac1{\lambda^2} g^e(\ona\wps^\infty, N)^2 + \frac1{\lambda^2} |\wka^\infty|^2 + 2\ \frac1{\lambda^2} |g^e(\ona\wps^\infty, N)|\ |\wka^\infty| \le 2 \ C (1+\sqrt{b})^2 \lb{sbw4}
 \end{align}
 Moreover, it follows from \eqref{pblowup} and \eqref{sbw1} that the points giving rise to the blow-up remain at finite distance from $0$.
 
 Let us observe also that \eqref{sbw1} and \eqref{sbw1a} give
\begin{align} &\wps(\wF(\cdot,\wtau))= \wps^\infty(\wF^\infty(\cdot(\tau(\wtau))) = \ln(r(\wF^\infty(\cdot,\tau(\wtau))^b) \lb{sbw1b}
\end{align}
that is, the function induced on $M$ by the immersions $\wF:M \flecha \re^2$ and $\wF^\infty: M \flecha \re^2$ is the same, let us call it $\psi_M$, but the metrics induced satisfy $\wg = \lambda^2 \wg^\infty$, which gives for the gradients of the above functions in the two different metrics the relation $\wna \psi_M = \lambda^{-2} \wna^\infty \psi_M$ (because $\wg^\infty(\wna^\infty \psi_M, X)= d \psi_M (X) =\wg(\wna \psi_M, X) = \lambda^2 \wg^\infty(\wna \psi_M, X)$). Moreover, the ordinary laplacians in these two metrics are related by $\wDe= \lambda^{-2} \wDe^\infty$.

From both expressions we obtain
\begin{align}
\wDe_{\psi_M} = \lambda^{-2} \wDe^\infty_{\psi_M}
\end{align}
Moreover:
\begin{align}
(\partial_\ws^m\wka_\psi)^2 = \frac1{\lambda^{2m+2}}(\partial_{\ws^\infty}^m\wka^\infty_\psi)^2, \qquad \parcial{\tau}{\wtau} = \frac1{\lambda^2},
\end{align}
and
\begin{align}
\ona^m \wps = \lambda^{-m} \ona^m\wps^\infty
\end{align}
From \eqref{E_igualdadprincipal} and the above expressions, we have 

\begin{align}
\parcial{}{\wtau}(\partial_{\ws}^m \wka_{\psi})^{2} &= \frac1{\lambda^2}\parcial{}{\tau}\( \frac1{\lambda^{2m+2}}(\partial_{\ws^\infty}^m \wka^\infty_{\psi})^{2}\) \nn \\
&= -\frac{2m+2}{\lambda^{2m+2}} (\partial_{\ws^\infty}^m \wka^\infty_{\psi})^{2}+ \frac1{\lambda^{2m+4}} \(\wDe^\infty_{\psi_M} \(\partial_{\ws^\infty}^m \wka^\infty_{\psi}\)^2
- 2  \(\partial_{\ws^\infty}^{m+1} \wka^\infty_{\psi}\)^2 \right. \nn \\
& \qquad \left. + 2  (a_{m0}  + a_{m1} \wka^\infty_{\psi} +a_{m2} (\wka^\infty_{\psi})^2) (\partial_{\ws^\infty}^m \wka^\infty_{\psi})^2   + \sum a_{iJ}   \((\wka_\psi^\infty)^i \partial_{\ws^\infty}^J \wka^\infty_{\psi}\) (\partial_{\ws^\infty}^m \wka^\infty_{\psi}) \)\nn \\
&= - (2m+2) (\partial_{\ws}^m \wka_{\psi})^{2}+  \wDe_{\psi_M} \(\partial_{\ws}^m \wka_{\psi}\)^2
- 2  \(\partial_{\ws}^{m+1} \wka_{\psi}\)^2  \nn \\
& \qquad  + 2  (a_{m0}  + a_{m1} \wka_{\psi} +a_{m2} (\wka_{\psi})^2) (\partial_{\ws}^m \wka_{\psi})^2   + \sum a_{iJ}   \((\wka_\psi)^i \partial_{\ws}^J \wka_{\psi}\) (\partial_{\ws}^m \wka_{\psi}) ,
\label{E_igualdadpb}
\end{align}
This equation is similar to that which appears in ordinary MCF for the blow-up in type I singularities (for instance, see \cite{man} page 59) then, standard arguments (see the same reference) show that 
\begin{align}\lb{numerito}
(\partial_{\ws}^m \wka_{\psi})^{2}\le D_m^2 \ \text{ for some constant $D_m$.} 
\end{align}

Because we are in the Euclidean plane and $r$ is the distance to a line, one has $|\ona r|=1$ and $\ona^mr=0$ for $m\ge 2$. From this and \eqref{sbw1a} it follows that $\ona^m \wps =  (-1)^{m+1} \fracc{(m-1)! b}{r^m} \ona r \otimes \overset{\overset{m}{\smile}}{\dots} \otimes \ona r$ and $|\ona^m \wps| \le  \fracc{(m-1)!b}{r^m}$. Moreover, from \eqref{1/r} it follows that the rescaled flows $F^j$	satisfy $\fracc1{r(F^j)}  \le \fracc1{\lambda_j}\frac{1}{\sqrt{C_4 (T-t_j-\lambda_j^{-2} \tau)}} = \frac{1}{\sqrt{C_4}} \frac{1}{\sqrt{C-\tau}}= \fracc{\sqrt{2}}{\sqrt{C_4}}\lambda(\tau)$. Then $\fracc{1}{r(\wF^\infty)} \le \fracc{\sqrt{2}}{\sqrt{C_4}}\lambda(\tau)$ and $\fracc{1}{r(\wF)} \le \fracc{\sqrt{2}}{\sqrt{C_4}}=: C_5$. This gives $|\ona^m \wps| \le  \ (m-1)!b\ C_5^m$. Then, writing  $\partial_{\ws}^m \<\ona\wps,N\>$ in function of the $\ona^\ell \wps$ and $\partial_s^J\wka_\psi$, we obtain from the above estimates that  $|\partial_{\ws}^m \wka| \le |\partial_{\ws}^m \wka_{\psi}| +|\partial_{\ws}^m \<\ona\wps,N\>|$ is bounded.

 Then there is a sequence of times $\wtau_n$ such that the $\wF(\cdot,\wtau_n)$ converges smoothly to a curve $\wF_\infty^\infty: M_\infty^\infty \flecha \re^2$.  
 
Now, let us check that we can apply formula \eqref{formon}  to the flow $\wF^\infty(\cdot,\tau)$. That is, we want to see that  $\partial (\wF^\infty(M_\infty,\tau)\cap C_R)$ is bounded by a finite natural number independent of $\tau$ and $R$ (for $\tau$ big enough). Let us consider, in the equivalent flow \eqref{egmcf}, the family $\wwr^j(z,\tau)$ with fixed $\tau$ which converges $C^\infty$ to some $\wwr^\infty(z,\tau)$. By Lemma \ref{discrt} there is a $t_0$ such that for every $t\ge t_0$  there exists a finite ordered family $z_1^j < z_2^j < ... < z_m^j$ of zeros of $\dot \wwr^j(\cdot, \tau)= \dot r(z/\lambda_j, t_j + \lambda_j^{-2} \tau)$. If we consider $\dot \wwr^t(\cdot, \tau)= \dot r(z/\lambda(t), t + \lambda(t)^{-2} \tau)$ as a function of $t\in[t_0,T[$, for every $\tau$ fixed, the corresponding zeros $z_1(t) < z_2(t) < ... < z_m(t)$ are continuous functions of $t$. Then, an argument similar to that given in the proof of Lemma 5.1 in \cite{aag} shows that $\lim_{j\to\infty} z^j_k = \lim_{t\to T}z_k(t)=:z_k$ exists for $k=1, ..., m$, 
where may be some $z_i$ are $-\infty$ or $\infty$.

Since $z_1^{j} < z_2^{j} < ... < z_m^{j}$, one has that $z_1 \le z_2 \le ... \le z_m$. Let $i_0= \min\{i;\ z_i> -\infty\}$ and $i_m= \max\{i;\ z_i < \infty\}$. Again by Lemma \ref{discrt}, the functions $\wwr^{j}$ are strictly monotone on the intervals $]-\infty, z_1^{j}[,\ ]z_1^{j},z_2^{j}[ ,\  ... ,\  ]z_m^{j},\infty[$, then their limit $\wwr^{\infty}$ is monotone on the intervals $]-\infty, z_{i_0}[,\ ]z_{i_0},z_{i_0+1}[ ,\  ... ,\  ]z_{i_m},\infty[$. Then, on each one of these intervals, the intersection of a line $r=constant$ with the graphics of $\wwr^\infty(\cdot,\tau)$ is one point or one segment. Then, the intersection of a line $r=constant$ with the graphics of $\wwr^\infty(\cdot,\tau)$ along all the domain of $\wwr^\infty$ is a finite union of segments and points. This implies that the boundary of the intersection of the graph of $\wwr^\infty(\cdot,\tau)$ with a square  centered at the origin consists on a finite number of points less than $2m$, and the condition in order formula \eqref{formon} be true when the manifold is not compact is fulfilled.

By formulae \eqref{formon}, \eqref{sbw1}, \eqref{sbw2a} and \eqref{sbw2aa}, taking into account that   $\widetilde\mu_\wtau = (2(T-t))^{-n/2} \mu^\infty_\tau$  and $\deri{}{\wtau} = \deri{\tau}{\wtau}\deri{}{\tau} = e^{-2 \wtau}  \deri{}{\tau} =\fracc1{\lambda^2}\deri{}{\tau} $,
\begin{align}
\deri{}{\wtau} \int_{M_\infty} e^{-|\wF|^2/2}  r(\wF)^b  \widetilde\mu_\wtau  &= \frac1{\lambda^2} \deri{}{\tau} \int_{M_\infty} e^{-|\wF^\infty|^2\lambda^2/2}  \lambda^b r(\wF^\infty)^b  \lambda^n  \mu^\infty_\tau = \frac{(2\pi)^{(n+b)/2}}{\lambda^2} \deri{}{\tau} \int_{M_\infty}\uu(\wF^\infty,\tau) r(\wF^\infty)^b\mu^\infty_\tau\nn \\
&=-\frac{(2\pi)^{(n+b)/2}}{\lambda^2}\int_{M_\infty}  (\kappa_\psi^\infty + \lambda^2\<\wF^\infty,N\> )^2  (\sqrt{2\pi}\lambda)^{-(n+b)} e^{-|\wF^\infty|^2 \lambda^2/2} r(\wF^\infty)^b \mu^\infty_\tau  \nn\\
&= - \frac{(2\pi)^{(n+b)/2}}{\lambda^2}\int_{M_\infty}   (\lambda \wka_\psi + \lambda \<\wF,\wN\> )^2  (\sqrt{2\pi})^{-(n+b)} \lambda^{n+b} e^{-|\wF/\lambda|^2\lambda^2/2} \lambda^{-b} r(\wF)^b \lambda^{-n} \widetilde\mu_\wtau \nn \\
&= - \int_{M_\infty}   ( \wka_\psi +  \<\wF,\wN\> )^2   e^{-|\wF|^2/2}  r(\wF)^b  \widetilde\mu_\wtau 
\end{align}

Then, for every $\wtau >0$,
\begin{align*}
-\int_\wtau^\infty \int_{M_\infty}  (\wka_\psi + \<\wF,N\> )^2  e^{-|\wF|^2/2} \widetilde \mu_{\psi \wtau} d\wtau & =  \int_\wtau^\infty  \parcial{}{\wtau}\int_{M_\infty} e^{-|\wF|^2/2} \widetilde \mu_{\psi \wtau} d\wtau \\
&= -\int_{M_\infty} e^{-|\wF|^2/2} \widetilde \mu_{\psi \wtau} + \lim_{\wtau\to\infty}\int_{M_\infty} e^{-|\wF|^2/2} \widetilde \mu_{\psi \wtau} .
\end{align*}
which is finite,
then $ \int_{M_\infty}  (\wka_\psi + \<\wF,N\> )^2  \widetilde \mu_{\psi \wtau} \us{\wtau\to\infty}{\flecha} 0$, and $M_\infty^\infty$ satisfies
$ \<\wF_\infty^\infty, N\>+ \wka_\psi = 0. $

When $b=m$, $\wka_\psi$ coincides with the mean curvature of a revolution hypersurface $M_\infty^\infty\times_{m \ln r} S^m$ of $\re^{m+2}$ and it is known by the classification of the mean convex shrinkers that it must be a Cylinder, then $M_\infty^\infty$ must be a line.
\end{demo}

 \bibliographystyle{alpha}

\noindent
Department of Mathematics \\
University of Valencia \\
46100-Burjassot (Valencia), Spain\\
 {miquel@uv.es} and {Francisco.Vinado@uv.es}

\end{document}